\input amstex
\documentstyle{amsppt}

\def\Spec{\roman{Spec}}
\define\sumk#1{{
        \sum_{\beta \colon k_r (\beta) = #1}
                \langle \Delta_2^{#1-1} \rangle_\beta   }}

\def\suma{
        \sum_{\beta \colon k_r (\beta) = 1}
                \langle \emptyset \rangle_\beta   }
\def\sumb{
        \sum_{\beta \colon k_r (\beta) = 2}
                \langle \Delta_2 \rangle_\beta   }

\def\adkr{k_r\circ}
\def\C{\bold{C}}

\def\P{\bold{P}}
\def\Q{\bold{Q}}

\def\Z{\bold{Z}}


\document
\magnification=1200
\NoBlackBoxes
\nologo
\vsize18cm



\bigskip

\centerline{\bf (SEMI)SIMPLE EXERCISES IN QUANTUM COHOMOLOGY }

\medskip

\centerline{\bf Arend Bayer${}^1$, Yuri I. Manin${}^2$}

\medskip

\centerline{\it ${}^1$Bonn University, Bonn, Germany}

\smallskip

\centerline{\it ${}^2$Max--Planck--Institut f\"ur Mathematik, Bonn, Germany}



\bigskip

{\bf Abstract.} The paper is dedicated to the study of
algebraic manifolds whose quantum cohomology or
a part of it is a semisimple Frobenius manifold.
Theorem  1.8.1 says, roughly speaking, that
the sum of $(p,p)$--cohomology spaces is a maximal
Frobenius submanifold that has chances to be semisimple.
Theorem 1.8.3 provides a version of the Reconstruction theorem,
assuming semisimplicity but not $H^2$--generation.
Theorem 3.6.1 establishes the semisimplicity for all
del Pezzo surfaces, providing an evidence for the conjecture
that semisimplicity is related to the existence
of a full system of exceptional sheaves of the
appropriate length. Finally, in \S 2 we calculate special
coordinates for three families of Fano threefolds with
minimal cohomology.

\bigskip

\centerline{\bf \S 0. Introduction}

\medskip

{\bf 0.1. Role of semisimplicity.} This paper is
a contribution to the study of algebraic
manifolds whose quantum cohomology (or at least its appropriate
subspace) is a generically semisimple
Frobenius manifold. This class of manifolds
is important in many respects. For example,
one expects that the higher genus correlators
and the correlators with gravitational descendants
for such manifolds can be calculated entirely
in terms of the Frobenius structure, when it is
semisimple (see [DuZh], [G2],
and the references therein). Unfolding spaces
of the hypersurface singularities (Landau--Ginzburg models) are generically
semisimple as well, therefore any mirror isomorphism
between such an unfolding space and (a part of)
quantum cohomology must involve semisimple
cohomology. 
Finally, to establish such an isomorphism,
we have to compute for both Frobenius manifolds
only finitely many constants (monodromy data, or
special coordinates) whereas one generally
needs infinitely many numbers in order
to describe a non--semisimple manifold (like
the quantum cohomology of a quintic).

\smallskip

In this paper, we focus
on the values of Dubrovin's canonical coordinates $(u_i)$
at the points of $H^2$ that is, parameter
space of the small quantum cohomology.
As soon as their identification
at a tame semisimple points is made, the remaining
part of the mirror picture requires a choice
of the flat metric compatible with multiplication. Unless it is explicitly
given, we have to recourse to one of the following
equivalent data describing it in an implicit way:

\smallskip

(i) Values of the diagonal coefficients of the flat metric
$\sum_i\eta_i (du^i)^2$ at a tame semisimple point
completed by the values of its first derivatives $\eta_{ij},$
that is, initial data for the 
second structure connection (cf. [Ma2], II.3).

\smallskip

(ii) Monodromy data for the first structure connection
and oscillating integrals for the deformed flat coordinates
(cf. [G1], [Du] and the references therein).

\smallskip
(iii) Choice of one of K.~Saito's primitive forms.

\smallskip
(iv) Choice of a filtration on the
cohomology space of the Milnor fiber (M.~Saito, cf. [He2]
and the references therein).

\smallskip
(v) Use of the semi--infinite Hodge structure 
(this is a far--reaching refinement of (iii), cf. [Bar]).

\smallskip

In this paper we discuss only (i).

\medskip

{\bf 0.2. General notation for quantum cohomology.} Let $V$ be a smooth projective manifold.
We choose a homogeneous basis of its cohomology space
$(\Delta_a),\,\Delta_a\in H^{|\Delta_a|}(V)$ where $a$
runs over a set of indices such that $\Delta_0$ always denotes
the dual class of $[V]$. The dual coordinates
on $H^*(V)$ are denoted $(x_a).$ The Poincar\'e form
$g_{ab}=(\Delta_a,\Delta_b)$ allows us to raise indices,
e.~g. $\Delta^a :=\sum_bg^{ab}\Delta_b.$  The letter $\beta$ generally denotes a variable
algebraic element of the group $H_2(V,\bold{Z})/(tors)$, and $q^{\beta}$ is the respective element
of the Novikov ring. We write $k$ or $k(V)$ for $c_1(V)$
and $k(\beta )$ for $(c_1(V),\beta ).$

\smallskip

The genus zero Gromov--Witten invariants are written
as correlators, that is polylinear functions indexed by $\beta$:
$$
H^*(V,\bold{Q})^{\otimes n}\to \bold{Q}:\ \Delta_{a_1}\otimes \dots\otimes
\Delta_{a_n}\mapsto \langle \Delta_{a_1} \dots 
\Delta_{a_n}\rangle_{\beta}.
$$
For $n=0$ we write the respective correlator as $\langle\emptyset\rangle_{\beta}.$ If $\roman{dim}\,V\ge 3,$
we have $\langle\emptyset\rangle_{\beta}=0.$

\smallskip

Sometimes it is convenient
to consider the total correlators with values in the completed Novikov ring:
$\langle \Delta_{a_1}, \dots ,
\Delta_{a_n}\rangle :=\sum_{\beta}\langle \Delta_{a_1}, \dots ,
\Delta_{a_n}\rangle_{\beta}\,q^{\beta}$
and extend them to the completed tensor algebra of $H^*(V)$.
The potential of quantum cohomology is the 
formal series in $(x_a)$ over the Novikov ring
which can be compactly written in this notation as $\Phi:=\langle e^{\sum_ax_a\Delta_a}\rangle$.
The quantum multiplication table is given by
$\Delta_a\circ\Delta_b=\sum_c\Phi_{abc}\Delta^c$ where
$\Phi_{abc}=\partial_a\partial_b\partial_c\Phi,\,
\partial_a=\partial/\partial x_a$. Using the standard
properties of the Gromov--Witten invariants (symmetry,
divisor axiom and identity axiom) we can
rewrite this in the following form convenient for
further calculations:
$$
\Delta_a\circ\Delta_b = \Delta_a\,\cup\,\Delta_b +
$$
$$
\sum_{\beta \ne 0}\sum_{c\ne 0}
\langle\Delta_a\Delta_b\Delta_c\,e^{\sum_{k:\,|\Delta_k|>2}x_k\Delta_k}
\rangle_{\beta}\Delta^c\,e^{(\sum_{k:\,|\Delta_k|=2}x_k\Delta_k,\beta )}q^{\beta}.
\eqno(0.1)
$$
Since $\beta\mapsto e^{(\sum_{k:\,|\Delta_k|=2}x_k\Delta_k,\beta )}q^{\beta}$
is a generic character of $H_2(V,\bold{Z})/(tors)$, as well as
$\beta\mapsto q^{\beta}$, we will allow ourselves in the future to
use (0.1) omitting $e^{(\sum_{k:\,|\Delta_k|=2}x_k\Delta_k,\beta )}$
and instead interpreting $q^{\beta}$ as a monomial function on the 
algebraic torus $T_V$ with the character group $H^2(V,\bold{Z})/(tors)$. 

\medskip

{\bf 0.3. Plan of the paper.} We start \S 1 with a review
of general properties of semisimple Frobenius manifolds.
Actually, a natural context for this discussion
is the realm of $F$--manifolds, that is,
''Frobenius manifolds without metric'', introduced in [HeMa].
We recall their definition and the main properties
related to the semisimplicity in 1.1--1.7.

\smallskip

In 1.8 we then prove the first main result of this paper,
Theorem 1.8.1.
Roughly, it says that the most natural candidate
for semisimple quantum cohomology is not the whole
cohomology of a Fano manifold, but only its $(p,p)$--part.
The second result, Theorem 1.8.3, is a new version
of the First Reconstruction Theorem in [KoMa].
Finally, we discuss the special coordinates for
Fano manifolds with minimal $(p,p)$--cohomology
(rank one in every dimension).

\smallskip

In \S 2 we provide explicit canonical coordinates
and the initial data for the second structure connection
for three families of Fano threefolds with minimal
$(p,p)$--cohomology. They can be useful
for finding Landau--Ginzburg superpotentials
and primitive forms for them. Since these threefolds are not
toric, no clear--cut prescription is available
for doing this.

\smallskip

The whole \S 3 is dedicated to the quantum cohomology
of the del Pezzo surfaces. This subject
was already treated in several papers, see [GP]
and the references therein. We systematically use the 
symmetry properties of the quantum cohomology with respect
to the Weyl groups, established in [GP]. This allows us to encode
a large amount of numerical data in rather compact
tables and, more important, to perform calculations
in quantum cohomology
in terms of the geometry of root systems and reflections.
Since the expected mirror picture involving the
Lagrangian vanishing cycles is also naturally formulated in this 
way, we feel that this is the right language.

\smallskip

The principal result of \S 3 is the Theorem 3.6,
establishing generic semisimplicity of the quantum
cohomology of all del Pezzo surfaces $V_r.$
For $r\le 3$ when $V_r$ are toric, the point
$q=1$ turns out to be tame semisimple.
This ceases to be true for $r\ge 5$, and instead we study
the behavior of quantum cohomology near
the boundary of the partial compactification
of the torus $T_{V_r}$, where the exponentiated coordinate
corresponding to an exceptional curve vanishes.

\smallskip

For this proof, the calculations of 3.3--3.5 are not needed,
and the reader may prefer to skip them. On the other hand,
they might be useful for finding the Landau--Ginzburg 
superpotentials in the non--toric cases.
\medskip

{\bf Acknowledgement.} We are thankful to Lothar
G\"ottsche who provided us with several values
of Gromov--Witten invariants missing in [GP], and to A.~Bondal,
who allowed us to use the tables 
of the Gromov--Witten invariants of
certain Fano threefolds, made by him, A.~Kuznetsov, and D.~Orlov.

\newpage

\centerline{\bf \S 1. Semisimplicity and $F$--manifolds}

\medskip 

\proclaim{\quad 1.1. Definition} An $F$--manifold
is a triple $(M,\circ ,e),$ where $M$
is a manifold and $\circ$ is an associative
supercommutative $\Cal{O}_M$--bilinear multiplication
$\Cal{T}_M\times \Cal{T}_M\to \Cal{T}_M$ with identity
vector field $e$, satisfying the
following axiom: for any (local) vector fields
$X,Y,Z,W$ we have
$$
[X\circ Y,Z\circ W]-[X\circ Y,Z]\circ W-
Z\circ [X\circ Y,W]
$$
$$
-X\circ [Y,Z\circ W] + X\circ [Y,Z]\circ W
 +X\circ Z\circ [Y,W]
$$
$$
-Y\circ [X,Z\circ W] + Y\circ [X,Z]\circ W +
Y\circ Z\circ [X,W] =0
\eqno(1.1)
$$
\endproclaim

\smallskip

This notion makes sense in any of the standard
categories of manifolds, and also supermanifolds,
if appropriate signs are introduced in the structural
identity, cf. [HeMa]. Below we will assume that $M$ is
a complex manifold.

\smallskip

A shorter formulation of (1.1) is
$$
\roman{Lie}_{X\circ Y}(\circ )=X\circ \roman{Lie}_{Y}(\circ )+
Y\circ \roman{Lie}_{X}(\circ )
\eqno(1.2)
$$
where $X,Y$ are any two local vector fields, and $\circ$
is considered as a tensor.

\medskip

{\bf 1.1.1. Question.} {\it Does the quantum multiplication
in $K$--theory (whose construction
is sketched by Givental in [G3]) satisfy (1.1)--(1.2)?}

\medskip

Returning to $F$--manifolds, we will say that 
$x\in M$ is a {\it semisimple point}
of the $F$--manifold $M$, if the following
equivalent conditions are satisfied:

\medskip

(i) {\it $T_xM$ endowed with the multiplication $\circ$
is a (commutative) semisimple algebra, that is,
isomorphic to $\bold{C}^n$, $n=\roman{dim}\,M,$ with
componentwise multiplication.}

\smallskip

(ii) {\it In a neighborhood of $x$, there exists a
basis of vector fields $e_i,\, i=1,\dots ,n,$ such that
$e_i\circ e_j=\delta_{ij} e_i$ and $[e_i,e_j]=0$ for all $i,j$.}

\medskip

$(M,\circ ,e)$ is called {\it generically semisimple},
if its set of semisimple points is open and dense.

\smallskip

>From (ii) it follows that in a neighborhood of any semisimple
point, $M$ admits a system of Dubrovin's canonical coordinates
$(u_1,\dots ,u_n).$
They have the following property: $e_i:=\partial /\partial u_i$
form a complete system of pairwise orthogonal idempotents
with respect to $\circ$. Any two local systems of canonical coordinates
differ by a constant shift $(u_j+c_j)$ followed by a permutation.

\smallskip

$F$--manifolds admit a natural product operation:
$$
(M_1,\circ_1,e_1)\times (M_2,\circ_2,e_2):=
(M_1\times M_2,\circ_1\boxplus \circ_2,,e_1\boxplus e_2) .
$$
The main structural identity (1.1) can be used as an 
integrability condition, in order to show the following
result ([He1], Theorem 4.2): {\it for any $x\in M$, 
the canonical decomposition
of $(T_xM,\circ )$ into a product of local algebras 
uniquely extends to the direct decomposition
of the germ of $(M,x)$ into the product of irreducible
germs of $F$--manifolds.}

\medskip

{\bf 1.2. Euler fields.} A semisimple point $x$ endowed with local canonical coordinates
$u_i$ is called {\it tame semisimple}, if $u_i(x)\ne u_j(x)$
for all $i\ne j.$
There are two basic setups in 
which the eventual constant shifts 
of canonical coordinates can be eliminated and the notion
of tame semisimplicity made coordinate--independent.

\smallskip

First, it might happen that the fundamental group of the
submanifold of semisimple points acts transitively
on the elements of a local idempotent system $(e_i)$.
Then $(u_i)$ can be normalized up to a common shift
$(u_i+c)$ and a permutation, and tameness or otherwise 
with respect to such
coordinates is defined unambiguously.

\smallskip

Second, it might happen that the $F$--manifold $(M,\circ ,e)$
is endowed with a structure Euler vector field $E$ (this is the case
of quantum cohomology). In the context of $F$--manifolds,
a local vector field $E$ is called an {\it Euler
field of (constant) weight $d_0$} if $\roman{Lie}_E(\circ )=d_0\circ$ that is,
for all local vector fields $X,Y$ we have
$$
[E,X\circ Y]-[E,X]\circ Y-X\circ [E,Y]=d_0X\circ Y.
\eqno(1.3)
$$
Euler fields form a sheaf of Lie algebras,
and weight is a linear function on this sheaf. The commutator
of two Euler fields has weight zero.

\smallskip

It is easy to describe all Euler fields in a neighborhood
of a semisimple point endowed with canonical coordinates
$(u_i).$ Namely, writing (1.3) for $X=e_i,\,Y=e_j,$ we see that any
local Euler field of weight $d_0$ is of the form
$$
E=d_0\sum_i u_ie_i +\sum_jc_je_j,
\eqno(1.4)
$$
where $c_j$ are arbitrary constants. Reversing this argument,
we see that given an Euler field $E$ of non--zero weight,
we can define uniquely up to permutation canonical coordinates
$(u_i)$ by the condition $E=d_0\sum_i u_ie_i$. Moreover, functions
$(d_0u_i)$ can be then invariantly described as the spectrum
of the operator $E\circ$ acting upon local vector fields.
Notice that in quantum cohomology $E$ is normalized in 
such a way that $d_0=1$, and we will adopt this convention.

\medskip

{\bf 1.3. The spectral cover.} Let $(M,\circ ,e)$ be an $F$--manifold.
By definition, its {\it spectral cover} is the relative
analytic spectrum $L:=\roman{Specan}\,(\Cal{T}_M,\circ )$, together
with two structure maps $\pi:\, L\to M$ and
$\sigma :\,\Cal{T}_M\to \pi_*(\Cal{O}_L)$. Clearly, $\sigma$ is an isomorphism
of sheaves of $\Cal{O}_M$--algebras. A point
$x\in M$ is semisimple iff $\pi$ is \'etale at $x$. 

\smallskip

Moreover, $L$ is endowed with a canonical closed embedding 
$i:\,L\to T^*M$ to the complex analytic
cotangent bundle of $M$. This embedding is induced by the
map of sheaves of algebras $S(\Cal{T}_M)\to (\Cal{T}_M,\circ )$
identical on $\Cal{T}_M$. 

\smallskip

The definition of $L$, as well as the definition
of a semisimple point, does not require $\circ$ to satisfy (1.1).
Generic semisimplicity is equivalent to the fact that
$L$ is reduced (has no nilpotents in the structure sheaf).
If this condition is satisfied, (1.1) becomes equivalent 
to the condition that $L$ is Lagrangian.

\medskip

{\bf 1.4. K.~Saito's frameworks and Landau--Ginzburg models.}
A large class of generically semisimple $F$--manifolds is furnished
by the construction due to K.~Saito in the theory
of singularities and axiomatized in [Ma1], [Ma2] (for
the richer Frobenius structure)
under the name of {\it K.~Saito's frameworks}. Briefly,
let $p:\,N\to M$ be a submersion  of complex
manifolds, generally non--compact. Then $N$ carries
the complex of sheaves of relative holomorphic forms
$(\Omega^*(N/M),d_{N/M}).$ Let $\Psi$ be a closed
relative 1--form, and $C\subset N$  the
closed analytic subspace $\Psi =0$ of $N$.
We can define a map $s:\,\Cal{T}_M\to p_*(\Cal{O}_C)$
in the following way. Let $X$ be a local vector field  on $M$, 
defined in a sufficiently small open subset $U.$
Cover $p^{-1}(U)$ by small open subsets $U_j$
in such a way that on each $U_i$ there exists
a lift $X_j$ of $X$. Then the maps
$X\mapsto i_{X_j}(\Psi )$ on $U_j$ glue together
and produce a well defined map of sheaves
$s:\, \Cal{T}_M\to p_*(\Cal{O}_C)$. (In [Ma1], [Ma2],
only the case $\Psi =d_{N/M}F$ was considered,
which is not a restriction locally on $N$ but
overlooks more global effects).

\smallskip

Assume that $s$ is an isomorphism of $\Cal{O}_M$--modules,
in particular, $p$ restricted to $C$ is flat.
Then $s$ induces a multiplication $\circ$
on $\Cal{T}_M$ coming from $p_*(\Cal{O}_C)$
and endowed with the identity vector field 
$e=s^{-1}(1)$. If $p:\,C\to M$ is generically \'etale,
$(M,\circ ,e)$ is an $F$--manifold. Its spectral cover
$L$ is then canonically identified with $C$, and
$s$ with $\sigma$. When these conditions are satisfied,
we will call $(p:\,N\to M,\Psi )$ a K.~Saito's framework
(in the context of $F$--manifolds).

\smallskip

If an $F$--manifold $(M,\circ ,e)$ 
is initially defined by another construction,
e.~g. is the quantum cohomology of a manifold $V$,
then any isomorphism of it with the base of some
Saito's framework $(p:\,N\to M,\Psi )$ is called {\it a Landau--Ginzburg representation}
of it. When $\Psi=d_{N/M}F$ for some local
function $F$ on $M$, $F$ is called the respective 
{\it LG--superpotential} of $(M,\circ ,e)$.
Canonical coordinates on a simply connected
$U\subset M$ over which $C$ is \'etale, are precisely the
critical values of $F$, that is, restrictions
of $F$ on various connected components of $p^{-1}(U).$

\smallskip

If $V$ is a Fano manifold with generically 
semisimple quantum cohomology, any Landau--Ginzburg
representation of it is traditionally considered
as a {\it mirror construction.}

\smallskip

\medskip

{\bf 1.5. Example: unfolding of isolated hypersurface
singularities.} Let $f=f(z_1,\dots,z_m)$ be the germ of
an analytic function with an isolated singularity at zero.
Let $g_1,\dots ,g_n$ be germs of analytic functions at zero
whose classes modulo the Jacobian ideal $J:= (\partial f/\partial z_i)$
form a basis of the Milnor ring of $f$. Put
$F=F(z,t)=f(z) +\sum_{j=1}^n t_jg_j(z)$. Consider $F$ as a function
on the germ $N$ at zero of the analytic space with
coordinates $(z,t)$ endowed with the projection
$(z,t)\mapsto (t)$ to the germ $M$ of the $(t)$--space.
Then $(N\to M, p, d_{N/M}F)$ is a K.~Saito's framework.
It produces a generically semisimple $F$--structure on the
base of miniversal deformation $M$ whose
spectral cover $C$ is smooth. Moreover, this (germ
of) $F$--manifold is {\it irreducible}, that is, cannot be represented
as a product of $F$--manifolds, because its tangent space
at zero is a local algebra.

\smallskip

Conversely, any irreducible germ of a generically
semisimple $F$--manifold {\it with
smooth spectral cover} $L$ is isomorphic to the miniversal
unfolding base of an isolated hypersurface
singularity, and this singularity is unique up to stable right
equivalence, that is, adding squares of extra $z$--coordinates
and performing invertible coordinate
changes. This beautiful theorem was proved in [He1]
(see Theorem 16.6), using earlier results by V.~Arnold and A.~Givental.

\smallskip

It gives a very neat intrinsic description of
irreducible germs of $F$--manifolds allowing the Landau--Ginzburg
representations of the well studied type.

\smallskip

If we only know that an $F$--manifold $(M,\circ ,e)$
has a smooth spectral cover $L$, then Hertling's
local decomposition theorem together with the result above
produces a rather detailed local information about possible
$LG$--potentials, but neither proves its existence,
nor characterizes it uniquely.

\medskip

{\bf 1.6. Example: $LG$--superpotentials of the quantum cohomology
of projective spaces.} Quantum cohomology
of $H^*(\bold{P}^r)$ is semisimple at zero.
The canonical coordinates of this point defined as the
eigenvalues of the multiplication by the standard
Euler field $E$ are simply
$$
u_k=\zeta^k (r+1),\ \zeta = e^{\frac{2\pi i}{r+1}}.
\eqno(1.5) 
$$

\smallskip

The respective K.~Saito's framework which is traditionally
considered as the mirror of $\bold{P}^r$ is constructed
as follows (cf. [G1], [Bar]). Denote by
$M$ the germ of $\bold{C}^{r+1}$ at zero, with
local coordinates $t_0,\dots,t_r .$ Denote by $N$ the direct
product of $M$ and the torus
$z_0\dots z_r=1$, and by $p$ its projection on $M$. 
Put $f(z):=z_0+\dots +z_r.$ Finally, consider the function
$$
F(z;t):= f(z)+ \sum_{j=0}^r t_j\,f(z)^j.
\eqno(1.6)
$$
An easy calculation shows that at $t=0$, 
the critical points of $F(z;0)=f(z)$ are
$z_0=\dots =z_r=\zeta^k,\ k=0,\dots ,r+1$, with the
respective critical values
$\zeta^k (r+1)$. Therefore, there exists a unique
isomorphism of the $F$--manifold $(M,\circ ,e)$
produced by this framework with the germ of the quantum
cohomology of $H^*(\bold{P}^r)$ extending the
above identification of the canonical coordinates at zero.

\smallskip

The deformation of $f(z)$ explicitly described by (1.6)
is formally similar to the usual construction of the unfolding space of the
isolated singularity: in fact, the first $r+1$ powers of $f(z)$
generate a basis of the Milnor ring of $f$ (which is semisimple
rather than local this time). However, theoretical reasons
for choosing this deformation are less compelling
than in the classical case: we are dealing with a global
but non--compact deformation problem.

\smallskip

For a thorough discussion of a class of such problems, 
see [Sa] and the references therein.
  
\medskip

{\bf 1.7. Metric and flat coordinates.} Let again $(M,\circ , e)$
be an $F$--manifold with semisimple base point.
Choose also an Euler field $E$ whose eigenvalues
constitute a normalized system of canonical coordinates $u_i$ on $M$.
An additional structure which turns $M$ into a Frobenius
manifold is the choice of a (complex analytic) flat
metric $g$ such that $g(X\circ Y,Z)=g(X, Y\circ Z)$
for all local vector fields $X,Y,Z.$ Dubrovin
has shown that in the canonical coordinates such a metric
can be expressed as $\sum_i e_i\eta\, (du_i)^2$ where
$\eta$ is a local function defined up to adding a constant.
For quantum cohomology, $g$ is the Poincar\'e form,
and one can take for $\eta$
the flat coordinate dual to the cohomology class of a point.
For the latter, cf. [Ma2], Ch. I,
Proposition 3.5 c).

\smallskip

If the base point of $M$ is tame semisimple, then
the metric can be uniquely reconstructed from
the values of $\eta_i=e_i\eta$ and $\eta_{ij}=e_ie_j\eta$
at this point. For a proof, see [Ma2], I.3.

\smallskip

The essential (and difficult) problem in this case
consists in studying the global analytic properties
of flat coordinates near the boundary of the tame semisimple
domain. An appropriate language for such a study
is the language of monodromy data for the first and the second structure connections of $(M,\circ , g)$ for which we refer to
[Du], [Guz], and references therein (first connection and
the Stokes data), and [Ma2], [He2] (second connection, 
Schlesinger's equations).

\medskip

{\bf 1.8. Semisimple quantum cohomology.} We now turn to the main
subject of this section, $F$--manifolds and Frobenius
manifolds of quantum cohomology.

\smallskip

First of all, semisimple Frobenius manifolds cannot have
odd coordinates. Hence we should restrict ourselves from the start by the even--dimensional quantum cohomology. It is well known that it is
endowed with the induced $F$--manifold structure, Euler field,
flat invariant metric and flat identity.

\smallskip

\proclaim{\quad 1.8.1. Theorem} a) Let $V$ be a projective manifold
whose even--dimensional quantum cohomology is generically semisimple
as a formal Frobenius manifold over the Novikov ring.
Then
$$
h^{p,q}(V)=0\ \roman{for\ all}\ p\ne q,\ p+q\equiv 0\,\roman{mod}\,2.
\eqno(1.7)
$$

b) Generally, for an arbitrary projective
manifold $V$, the formal completion at of $\oplus_{p\in\bold{Z}} H^{p,p}(V)$
at the point of the classical limit
is endowed with an induced structure of Frobenius submanifold
with flat identity and Euler field.
\endproclaim

\smallskip

Recall that at the point of the classical limit
vanish all the flat coordinates outside $H^2$ and exponentiated
coordinates on $H^2$. More precisely, the $H^2$ component
of this  point lies
in the partial compactification of the respective torus
corresponding to the K\"ahler cone. 

\smallskip

{\bf 1.8.2. Dubrovin's conjecture and a generalization.} Classical calculations 
show that there exist Fano complete intersections in projective
spaces for which (1.7) fails. For example, $2m$--dimensional intersections
of three quadrics have 
$$
h^{m-1,m+1}=h^{m+1,m-1}=m(m+1)/2.
$$
Hence their
even--dimensional quantum cohomology cannot be semisimple.

\smallskip

Thus, the best result that one can expect is that $\oplus_{p\in\bold{Z}} H^{p,p}(V)$
is semisimple for an interesting class of Fano manifolds.  

\smallskip

B.~Dubrovin ([Du], 4.2.2) conjectured that the even quantum
cohomology of $V$ is semisimple {\it if and only if} $V$ is Fano 
and $Der^b(Coh\,V)$ admits
a full system of exceptional objects whose length equals
$\roman{dim}\,H^*(V)$. This agrees with Theorem 1.8.1,
because on such a manifold all cohomology
classes are algebraic, and in particular, $H^{odd}(V)=0$ and
$h^{p,q}=0$ for $p\ne q$.  Theorem 1.8.1 suggests the following
strengthening of this conjecture: {\it $\oplus_{p\in\bold{Z}} H^{p,p}(V)$
(with quantum multiplication) is semisimple if and only if $V$ is Fano 
and $Der^b(Coh\,V)$ admits
a full system of exceptional objects whose length equals $\sum_p h^{p,p}(V).$} 

\smallskip

For some  positive results on the semisimplicity in a more
narrow sense, see [Beau], [TX], and the references in these papers.
In particular, in [TX] it is shown that for certain
Fano complete intersections the operator of the quantum multiplication
by the canonical class is generically semisimple
on the classical subring of $\oplus_{p\in\bold{Z}} H^{p,p}(V)$
generated by $H^{1,1}(V)$. However, we do not know,
whether this subring is stable with respect to the quantum
multiplication. 

\medskip

{\bf Proof of the Theorem 1.8.1.} a) It is based on the consideration
of Euler fields, discussed in 1.2. Namely,
from (1.4) one easily deduces that if on a generically
semisimple $F$--manifold two Euler fields of non--zero weights
commute, they must be proportional.

\smallskip

Now let us turn to the quantum cohomology.
As Sheldon Katz remarked, it admits generally
at least two commuting Euler fields. In order to write them down explicitly,
choose a bigraded homogeneous basis $(\Delta_a )$ of $H=H^{*,*}(V,\bold{C})$
considered as the space of flat vector fields, and let
$(x_a)$ be the dual flat coordinates vanishing at zero. 
Let $\Delta_a\in H^{p_a,q_a}(V).$ 
Put $-K_V=\sum_{p_b+q_b=2}r_b\Delta_b.$ Then 
$$
E_1:=\sum_a (1-p_a)x_a\Delta_a+\sum_br_b\Delta_b,
\eqno(1.8)
$$
$$
E_2:=\sum_a (1-q_a)x_a\Delta_a+\sum_br_b\Delta_b
\eqno(1.9)
$$
are Euler of weight 1. To check this, one can repeat
the calculations made in [Ma2], III.5.4 with obvious changes. 
In the final count,
the reason for this is that Gromov--Witten correspondences
are algebraic, and therefore are represented by the cohomology
classes of types $(p,p).$

\smallskip

Assuming now that $h^{p,q}(V)\ne 0$ for some $p\ne q,\,p+q\equiv 0\,
\roman{mod}\,2,$ we see that the restrictions of (1.8) and (1.9)
to the even--dimensional cohomology cannot be proportional.
This contradiction concludes the proof of a).

\smallskip

b)  With the same notation, a correlator
$\langle \Delta_{a_1}\dots \Delta_{a_n}\rangle_{\beta}$
vanishes unless $\sum_i p_{a_i} = \sum_i q_{a_i}$ because the
Gromov--Witten correspondences are represented by algebraic cycles.
Let $\Phi$ be the potential of the quantum cohomology.
Assume that $p_a=q_a$ and $p_b=q_b$, but $p_c \ne q_c$
for some $a,b,c$. From the preceding remark it follows that
$\Phi_{ab}^c$ vanishes after restriction to the subspace
$\oplus_p H^{p,p}$ given by the equations $x_e=0$ for all
$e$ such that $p_e\ne q_e.$ Hence the $\circ$ product of two
vector fields tangent to this submanifold is again tangent to it.
It is also clear that the Poincar\'e form restricts
to a non--degenerate flat metric, and the identity field
and the Euler field are tangent to it as well.

\medskip

We will now prove a theorem which establishes
a result of the same type as the First Reconstruction
Theorem of [KoMa], but under the different assumptions.

\medskip

\proclaim{\quad 1.8.3. Reconstruction Theorem} Assume that 
the $(p,p)$--part of the quantum cohomology of $V$ is generically semisimple
and, moreover, admits a tame semisimple point
lying in the subspace $H^2(V)$ (parameter space of the small
quantum deformation).

\smallskip

In this case all genus zero Gromov--Witten invariants
of $\oplus_p H^{p,p}(V)$ can be reconstructed from the 
correlators $\langle \gamma_1\otimes\dots
\otimes\gamma_n\rangle_{\beta}$ with $n\le 4$, $\gamma_i\in H^{a_i,a_i}(V)$,
and $\beta \in H_2(V,Z)$ which can be non--vanishing
only if
$$
k(\beta ):=(c_1(V),\beta )=\sum_{i=1}^n (a_i-1) +3 -\roman{dim}\,V.
\eqno(1.10)
$$
In particular, if $\pm c_1(V)$ is numerically effective, finitely
many correlators suffice for the complete reconstruction.
\endproclaim

\smallskip

{\bf Proof.} As was explained in 1.7,
the whole germ of the Frobenius manifold
at such tame semisimple point of $H^2$ is determined by its 
canonical coordinates $(u_i^0)$,
the values $\eta_i^0$ at this point of the
coefficients $\eta_i$ of the metric (Poincar\'e form) $\sum \eta_i(du_i)^2$,
and the values $\eta_{ij}^0$ of their first derivatives $\eta_{ij}=e_i\eta_j.$

\smallskip

One easily sees that in order to calculate these data,
it suffices to know the expressions of $(u_i)$ through some
flat coordinates modulo $J^2$ where $J$ is the ideal of
equations of $H^2$ inside $\oplus H^{p,p}.$ Now, $(u_i)$
are the eigenvalues of the operator $E\circ$ acting on vector fields.
To calculate them modulo $J^2$, it suffices to know
the potential $\Phi$ of quantum cohomology modulo $J^5$
because the structure constants of the $\circ$--multiplication
are the third derivatives of the potential. With this 
precision, $\Phi$ is determined by all four--point correlators.
Of course, (1.10) is a special case of the Dimension Axiom.

\medskip

{\bf 1.9. Fano manifolds with minimal $(p,p)$--cohomology.} The finite
family of numbers $(u^0_i, \eta_{i}^0,\eta_{ij}^0)$ 
essentially coincides with what was called {\it special
coordinates} of the tame semisimple germ of Frobenius manifold,
cf. [Ma2], III.7.1.1. In this subsection, we will show
how to calculate them for the $\oplus H^{p,p}(V)$--part
of the quantum cohomology of those Fano manifolds, for which $\roman{dim}\,H^{p,p}(V)=1$
for all $1\le p \le \roman{dim}\,V:=r.$ This generalizes
our old computation for projective spaces.

\smallskip

We will work with a homogeneous basis  $\Delta_p\in H^{p,p}(V)$
consisting of rational cohomology classes
satisfying the following conditions:
$\Delta_0$ = the dual class of $[V]$, $\Delta_1 = c_1(V)/\rho$
is the ample generator of $\roman{Pic}\,V,$ and $\rho$
is called {\it the index} of $V.$ Furthermore,
$\Delta_{r-p}$ is dual to $\Delta_p$, that is $(\Delta_p,\Delta_{r-p})=1$,
$\Delta_r$ = the dual class of a point. The dual coordinates
are denoted $x_0,\dots ,x_r.$ The nonvanishing
$\beta =0$ correlators are coefficients of the
cubic self--intersection form
$$
(x_0\Delta_0 +\dots +x_r\Delta_r)^3.
$$
We put
$$
[d;\,a_1, \dots ,a_k]:=
\langle \Delta_{a_1}\dots \Delta_{a_k}\rangle_{d\Delta_{r-1}}\, .
$$
These symbols satisfy the following relations:

\medskip

(i) If $r\ge 3$, $[d;\,a_1,\dots ,a_k]\ne 0$ and $d>0$, then  necessarily
$$
k>0,\ a_i>0\ \roman{for\ all}\ i,\ \roman{and}\ d\rho =\sum_{i=1}^k (a_i-1)
+3-r 
\eqno(1.11)
$$
(see (1.10)). 

\smallskip

(ii) $[d;\,a_1,\dots ,a_k]$ is symmetric with respect to
the permutations of $a_1,\dots ,a_k .$

\smallskip

(iii) $[d;\,1, a_2,\dots ,a_k]=d\,[d;\,a_2,\dots ,a_k].$

\smallskip

(iv) Associativity relations, expressing the associativity of the
multiplication (0.1).

\medskip

The multiplication table in the first infinitesimal neighborhood
of $H^2(V)$ involves only up to four--point correlators and
looks as follows:
$$
\Delta_a\circ\Delta_b = \Delta_a\,\cup\,\Delta_b +
$$
$$
\sum_{d\ge 1}\sum_{c\ge 1} \left(
[d;\,a,b,c] + \sum_{f\ge 2} \,[d;\,a,b,c,f]\,x_f\right)\,
\Delta_{r-c} q^d\, .
\eqno(1.12)
$$
Finally, the (restricted) Euler field of weight 1 is
$$
E=\sum_{p=0}^r (1-p)\,x_p \Delta_p + \rho\,\Delta_1 
\eqno(1.13)
$$

\smallskip

As we already remarked, under the assumptions 
of  the Theorem 1.8.3, the eigenvalues of $E\circ$
at the generic point  of
$H^2$ are pairwise distinct and determine
the canonical coordinates of this point. 
We will have to calculate in the first infinitesimal neighborhood
of $H^{2}$ and therefore we will consider all the relevant
quantities as consisting of two summands: restriction
to $H^2$ and the linear (in $x_a$) correction term, in particular
$$
u_i:=u_i^{(0)}+u_i^{(1)} .
$$ 
The remaining special coordinates are given by the following formulas.
\medskip

\proclaim{\quad 1.9.1. Theorem} Put 
$$
e_i := \frac{\prod_{j\ne i} (E-u_j)}{\prod_{j\ne i} (u_i-u_j)}=
e_i^{(0)} + e_i^{(1)}\, .
\eqno(1.14)
$$
where the multiplication is understood in the sense
of quantum cohomology with the coefficient ring
extended by $(u_i)$ and $(u_i-u_j)^{-1}.$

\smallskip

Then we have on $H^2$:
$$
\eta_i = e_i^{(0)} (x_r),\ \eta_{ij} =e_i^{(0)} e_j\, (x_r)
\eqno(1.15)
$$
where in (1.15) $e_i$ are considered as vector fields
acting upon coordinates via $\Delta_a=\partial/\partial x_a.$
\endproclaim

\smallskip
  
{\bf Proof.} The elements $e_i$ are the basic pairwise orthogonal
idempotents in the quantum cohomology ring at the
considered point satisfying $E\circ e_i=u_ie_i.$
The metric potential $\eta$ is $x_r$. 

\smallskip

Here is an efficient way of computing $e_i^{(1)}.$
First, compute $\omega_{i}$ defined by the identity
in the first neighborhood:
$$
e_i^{(0)}\circ e_i^{(0)}=e_i^{(0)} +\omega_{i}.
\eqno(1.16)
$$
Then we have
$$
e_i^{(1)}=-\frac{\omega_i}{2e^{(0)}_i-1} = \omega_i \circ (1-2e^{(0)}_i)
\eqno(1.17)
$$
where the division resp. multiplication is again made in the first
neighborhood.

\smallskip

In fact, this follows from 
$$
(e_i^{(0)}+e_i^{(1)})^{\circ 2} =
e_i^{(0)}+e_i^{(1)}
$$
and (1.16).

\newpage

\centerline{\bf \S 2. Fano threefolds with minimal cohomology}

\medskip

{\bf 2.1. Notation.} Let $V$ be a Fano threefold. 
We keep the general notation of 1.3, but now
consider only the case $r=3$. Besides the index $\rho$,
we consider {the degree} $\delta :=(c_1(V)^3)/\rho^3$ of $V.$ 

\smallskip
 
There exist four families
of Fano threefolds $V=V_{\delta}$ with cohomology $H^{p,p}(V,\bold{Z})\cong\bold{Z}$
for $p=0,\dots ,3$ and $H^{p,q}(V,\bold{Z})=0$ for $p\ne q.$
Besides $V_1=\bold{P}^3$ and the quadric $V_2=Q$, they are
$V_5$ and $V_{22}$, with degree as subscript; their
indices are, respectively, 4,\,3,\,2,\,1. One can get
a $V_5$ by considering a generic codimension three linear section of 
the Grassmannian of lines in $\bold{P}^5$ embedded
in $\bold{P}^9.$ 

\smallskip

The nonvanishing
$\beta =0$ correlators are coefficients of the
cubic self--intersection form
$$
(x_0\Delta_0 +\dots +x_3\Delta_3)^3=\delta\,x_1^3+3x_0^2x_3+x_0x_1x_2\,.
$$

In this section, we will deal only with $Q$, $V_5$ and $V_{22}$, since
projective spaces of any dimension were treated by
various methods earlier: see [Ma3] for special coordinates,
[Du], 4.2.1 and [Guz] for monodromy data, [Bar]
for semiinfinite Hodge structures.

\smallskip

{\bf 2.2. Tables of correlators.} The following tables provide 
the coefficients of the multiplication table (1.12).

\smallskip

It suffices to tabulate the primitive correlators,
where primitivity means that
$a_i > 1$ and $a_i\le a_{i+1}.$ The symmetry and the
divisor identities furnish the remaining correlators.

\bigskip

{\it Manifold\ $Q$:}

\medskip

\settabs 8\columns
\+ [1;\,2,3] & [1;\,2,2,2] & [2;\,3,3,3]& [2;\,2,2,3,3] 
& & & & &\cr

\+ 1 & 1 & 1& 1 & & & & & \cr

\bigskip

{\it Manifold\ $V_5$:}

\medskip

\settabs 8\columns
\+ [1;\,3] & [1;\,2,2] & [2;\,3,3]& [2;\,2,2,3] 
&[3;\,3,3,3] & [2;\,2,2,2,2]& [3;\,2,2,3,3]&[4;\,3,3,3,3] &\cr

\+ 3 & 1 & 1& 1 & 1& 1& 2 & 3 &\cr

\bigskip

{\it Manifold\ $V_{22}$:}

\medskip

\settabs 8\columns
\+ [1;\,2] & [2;\,3] & [2;\,2,2]& [3;\,2,3] 
&[4;\,3,3] & [3;\,2,2,2]& [4;\,2,2,3]&[5;\,2,3,3] &\cr

\+ 2 & 6 & 1& 3 & 10& 1&  4& 16 &\cr

\medskip

\settabs 8\columns

\+ [6;\,3,3,3] & [4;\,2,2,2,2] & [5;\,2,2,2,3]& [6;\,2,2,3,3] 
&[7;\,2,3,3,3] & [8;\,3,3,3,3]& & &\cr

\+ 65 & 2 &9 & 41& 186 & 840& &  &\cr

\bigskip

The tables were compiled in the following way.
First, (1.10) furnishes the list of all primitive correlators
that might be (and actually are) non--vanishing.
Second, several correlators corresponding
to the smallest values of $n$ in (1.10) must be computed
geometrically: $n=2$ for $Q$, $n=1,2$ for $V_5$,
and $V_{22}$. These values were computed by A.~Bondal,
D.~Kuznetsov and D.~Orlov. Third, the associativity
equations uniquely determine all the remaining
correlators, in the spirit of the First Reconstruction
Theorem of [KoMa].

\medskip

{\bf 2.3. Canonical coordinates.} The canonical coordinates
on $H^2\,\cap\,\{x_0=0\}$ expressed in
terms of the flat coordinates are the roots
$u_0, \dots , u_3$ of the following characteristic equations
of the operator $E\circ$:
$$
Q:\qquad u^4-108\,q\,u=0,
\eqno(2.1)
$$
$$
V_5: \qquad u^4-44\,q\,u^2-16\,q^2=0,
\eqno(2.2)
$$
$$
V_{22}:\qquad (u+4\,q)\,(u^3-8\,qu^2-56\,q^2u-76\,q^3)=0.
\eqno(2.3)
$$
If $x_0\ne 0,$ one must simply add $x_0$ to the values above.

\medskip

{\bf 2.3.1. Question.} Find ``natural'' functions $f(z)$
whose critical values at $0$ are roots of (2.1)--(2.3)
and whose unfolding space carries an appropriate flat metric.

\bigskip

{\bf 2.4. Multiplication tables, idempotents, and metric
coefficients.} The remaining special coordinates $\eta_i,\,\eta_{jk}$
were calculated using the 
multiplication tables in the first neighborhood of $H^2$
obtained by specializing (1.12).
We calculated $e_i^{(0)}$ by determining the eigenvectors of $\text{ad}\ E$;
then we used equation (1.17) to get $e_i^{(1)}$.

\bigskip

{\it Manifold\ $Q$:}

\medskip

$$
\Delta_1^2=2\,\Delta_2+q\,\Delta_1x_3+q\,\Delta_0x_2,
$$
$$
\Delta_1\Delta_2=\Delta_3+q\,\Delta_0 +q\,\Delta_2x_3+q\,\Delta_1x_2,
$$ $$
\Delta_1\Delta_3 = q\,\Delta_1 + q\,\Delta_2x_2+2q^2\Delta_0x_3,
$$ $$
\Delta_2^2 = q\,\Delta_1 + q\,\Delta_2 x_2 + q^2\Delta_0x_3,
$$ $$
\Delta_2\Delta_3 = q\,\Delta_2 + q^2\Delta_0x_2 + q^2\Delta_1x_3,
$$ $$
\Delta_3^2 = q^2\Delta_0 + q^2\Delta_1x_2 + 2q^2\Delta_2x_3.
$$

\smallskip

Let $\xi_i$, $i=1,\dots ,3$ be the three roots of $\xi^3=4q$.
Then the respective idempotents have the following form:

$$
e_0 = \frac{1}{2}\Delta_0 - \frac{1}{2q} \Delta_3
	+ \frac{x_2}{4}\Delta_1 + \frac{x_3}{2}\Delta_2
$$ 
and, for $i = 1, \dots ,3$,
$$
e_i = \frac{1}{6}\Delta_0 + \frac{\xi_i^2}{12q} \Delta_1
	+ \frac{\xi_i}{6q} \Delta_2 + \frac{1}{6q} \Delta_3
- \frac{\xi_i x_2}{36}\Delta_0 
$$ 
$$
		- \left(\frac{x_2}{12} + \frac{\xi_i x_3}{12}\right) \Delta_1
	-\left(\frac{\xi_i^2 x_2}{18q} + \frac{x_3}{6}\right) \Delta_2
	-\left(\frac{\xi_i x_2}{27q} + \frac{\xi_i^2 x_3}{12q}\right) \Delta_3.
$$

Probably the most direct and efficient test of our computations is to
simply compute the pairwise products of these idempotents using the
multiplication table above. This verifies the formulas of the $e_j$ while
checking at the same time that our multiplication table yields an
associative product.

\smallskip

Note that $\Delta_i = \frac{\partial}{\partial x_i}$ and
$\frac{\partial}{\partial x_1} q = q$ since we identify $q$ with $e^{x_1}$
to get a proper (non--formal) Frobenius manifold. So we get as
special coordinates (where $i,j \in \{1,2,3\}$):

$$\eta_0 = - \frac{1}{2q}, \quad \eta_i = \frac{1}{6q}, \quad
\eta_{00} = 0,
$$ $$
\eta_{i0} = \frac{\xi_i^2}{12q}\Delta_1 \left(\frac{-1}{2q} \right)
	= \frac{\xi_i^2}{24q^2}, \quad
\eta_{0i} = \frac{-1}{2q} \Delta_3 \left(-\frac{\xi_i^2 x_3}{12q}\right)
	= \frac{\xi_i^2}{24q^2},
$$ $$
\eta_{ij} = \frac{\xi_i^2}{12q}\Delta_1 \left( \frac{1}{6q} \right)
	- \frac{\xi_i}{6q} \Delta_2 \left( \frac{\xi_j x_2}{27q} \right)
	- \frac{1}{6q} \Delta_3 \left(\frac{\xi_j^2 x_3}{12q} \right)
	= - \frac{\xi_i^2}{72q^2} - \frac{\xi_i \xi_j}{162q^2}
		- \frac{\xi_j^2}{72q^2}.
$$
The symmetry of $\eta_{ij}$ is just an additional check of our computations,
as we know in general that $e_i$ are commuting vector fields.

\bigskip

{\it Manifold\ $V_5$:}

\medskip

$$
\Delta_1^2 = 5 \Delta_2 + 3q\,\Delta_0 + 3q\,\Delta_2x_3 + q\,\Delta_1x_2
		+ 4q^2\Delta_0x_3,
$$ 
$$
\Delta_1\Delta_2 = \Delta_3 + q\,\Delta_1 + q\,\Delta_2x_2 + 2q^2\Delta_1x_3
			+ 2q^2\Delta_0x_2,
$$ 
$$
\Delta_1\Delta_3 = 3q\,\Delta_2 + 2q^2\Delta_0 + 4q^2\Delta_2x_3 
			+ 2q^2\Delta_1x_2 + 3q^3\Delta_0x_3,
$$ 
$$
\Delta_2^2 = q\,\Delta_2 + q^2\Delta_0 + 2q^2\Delta_2x_3
		+ q^2\Delta_1x_2 + 2q^3\Delta_0x_3,
$$ 
$$
\Delta_2\Delta_3 = q^2\Delta_1 + 2q^2\Delta_2x_2 + 2q^3\Delta_1x_3
			+2q^3\Delta_0x_2,
$$ 
$$
\Delta_3^2 = 2q^2\Delta_2 + q^3\Delta_0 + 3q^3\Delta_2x_3
		+2 q^3\Delta_1x_2 + 3 q^4\Delta_0x_3 .
$$

Let $u_i$ be the roots of $u^4 - 44\,q\,u^2 -16\,q^2$. The idempotents are
given by
$$
4000 q^3 e_i =
1440{q}^{3}\Delta_0
-20 q^2{u_i}^{2}\Delta_0
+70q u_i^{3}\Delta_1
-3040{q}^{2}u_i\Delta_1
$$ $$
-880{q}^{2}\Delta_{2}
+40q{u_i}^2\Delta_2
+4920q u_i\Delta_3
-110{u_i}^3 \Delta_3
$$ $$
-1968u_i{q}^{3}x_2\Delta_0
+44{u_i}^{3}{q}^{2}x_2\Delta_0
+352{q}^{4}x_3\Delta_0
-16 {u_i}^{2}{q}^{3}x_3\Delta_0
$$ $$
+176{q}^{3}x_2\Delta_1
-8{q}^{2}{u_i}^{2} \Delta_1 x_2
-5412u_i{q}^{3}\Delta_1x_3
+121{u_i}^{3}{q}^{2}\Delta_1x_3
$$ $$
-2864 u_i{q}^{2}x_2\Delta_2
+62 q{u_i}^{3}x_2\Delta_2
+1056{q}^{3}\Delta_2x_3
-48x_3{u_i}^{2}{q}^{2}\Delta_2
$$ $$
-16 q u_i^{2}x_2\Delta_3
+6036{q}^{2}u_i x_3\Delta_3
+352{q}^{2}x_2\Delta_3
-138q{u_i}^{3}x_3\Delta_3.
$$

The special coordinates $\eta_{ii}$ are
$$
\eta_{ii} 
= {\frac {-964\,q+21\,{u_{{i}}}^{2}}{800\,{q}^{3}}}
= \frac {-{251}\pm 105\sqrt{5}}{400\,q^2}.
$$

Now since the Galois group of $u^4 - 44\,q\,u^2 -16\,q^2$ obviously does
not act transitively on the pairs of roots, we have to distinguish two cases
in determining the $\eta_{ij}$. So we fix a root $u_1$; the other roots
are given by $u_2=-u_1$ and $u_{3,4}^2 - 11\,q + u_1^2 = 0$. We calculated

$$
\eta_{12} = \eta_{21} = 
{\frac {932\,q-21\,{u_{{1}}}^{2}}{800\,{q}^{3}}}
= {\frac {47\mp 21\,\sqrt {5}}{80\,{q}^{2}}} \quad \text{and}
$$
$$
\eta_{13} =
{\frac {- 3\,u_{{3}}u_{{1}} + 4\,q}{200{q}^{3}}}.
$$

All other coordinates are obtained from these via Galois permutations.

\bigskip

{\it Manifold\ $V_{22}$:}

\medskip

$$
\Delta_1^2 = 22\Delta_2 + 2q\,\Delta_1 + 24q^2\Delta_0
$$ $$
		+ 2q\,\Delta_2x_2 + 48q^2\Delta_2x_3 + 4q^2\Delta_1x_2
		+27q^3\Delta_1x_3 + 27q^3\Delta_0x_2 + 160q^4\Delta_0x_3,
$$ $$
\Delta_1\Delta_2 = \Delta_3 + 2q\,\Delta_2 + 2q^2\Delta_1 + 9q^3\Delta_0
$$ $$
		+ 4q^2\Delta_2x_2 + 27q^3\Delta_2x_3 + 3q^3\Delta_1x_2
		+ 16q^4\Delta_1x_3 + 16q^4\Delta_0x_2 + 80q^5\Delta_0x_3,
$$ $$
\Delta_1\Delta_3 = 24q^2\Delta_2 + 9q^3\Delta_1 + 40q^4\Delta_0
$$ $$
		+ 27q^3\Delta_2x_2 + 160q^4\Delta_2x_3 + 16q^4\Delta_1x_2
		+ 80q^5\Delta_1x_3 + 80q^5\Delta_0x_2 + 390q^6\Delta_0x_3,
$$ $$
\Delta_2^2 = 2q^2\Delta_2 + q^3\Delta_1 + 4q^4\Delta_0
$$ $$
		+ 3q^3\Delta_2x_2 + 16q^4\Delta_2x_3 + 2q^4\Delta_1x_2
		+ 9q^5\Delta_1x_3 + 9q^5\Delta_0x_2 + 41q^6\Delta_0x_3,
$$ $$
\Delta_2\Delta_3 = 9q^3\Delta_2 + 4q^4\Delta_1 + 16 q^5\Delta_0
$$ $$
		+ 16q^4\Delta_2x_2 + 80q^5\Delta_2x_3 + 9q^5\Delta_1x_2
		+ 41q^6\Delta_1x_3 + 41q^6\Delta_0x_2 + 186q^7\Delta_0x_3,
$$ $$
\Delta_3^2 = 40q^4\Delta_2 + 16q^5\Delta_1 + 65q^6\Delta_0
$$ $$
		+ 80q^5\Delta_2x_2 + 390q^6\Delta_2x_3 + 41q^6\Delta_1x_2
		+ 186q^7\Delta_1x_3 + 186q^7\Delta_0x_2 + 840q^8\Delta_0x_3.
$$

\smallskip

The first idempotent is given by
$$
e_0 = \frac{\Delta_0}{2} + \frac{2\Delta_2}{q^2} - \frac{\Delta_3}{2q^3}
+q x_2 \Delta_0+2\,{q}^{2}{x_3} \Delta_0
$$ $$
+\frac{x_2}{4}\,{\Delta_1}
+\frac{q x_3}{2}{\Delta_1}
+{\frac {2\,{ x_2}}{q}}{\Delta_2}
+\frac{3 x_3}{2}{ \Delta_2}
-{\frac{x_2}{2 q^2}}{\Delta_3}.
$$

Now let $u_i, i=1,\dots ,3$ be the roots of $u^3-8 q u^2 - 56 q^2 u - 76 q^3$.
Then the respective idempotents are given by

$$
5324\,q^4\,e_i = 
$$ $$
\left(
-71742{q}^{4}
-24552 {q}^{3}{u_i}
+2354{q}^{2}{{u_i}}^{2}\right)\Delta_0
+\left(
-30272{q}^{3}
-10186{q}^{2}{ u_i}
+979q{{u_i}}^{2}\right)\Delta_1
$$ $$
+\left(
-118712{q}^{2}
-38126q{u_i}
+3696{{u_i}}^{2}\right){\Delta_2}
+\left(
49346q
+16192{u_i}
-1562{\frac {{{u_i}}^{2}}{q}}\right){\Delta_3}
$$ $$
+\left(
-130876{q}^{5}
- 43283{u_i}{q}^{4}
+4168{{u_i}}^{2}{q}^{3}\right){x_2}\Delta_0
$$ $$
+\left(
- 483464{q}^{6}
-161648{u_i}{q}^{5}
+15528 {{u_i}}^{2}{q}^{4}\right){x_3}\Delta_0
$$ $$
+\left(
-38977{q}^{4}
-12940{u_i}{q}^{3}
+1245{{u_i}}^{2}{q}^{2}\right){x_2}{\Delta_1}
$$ $$
+\left(
-145898{q}^{5}
-48889{u_i}{ q}^{4}
+4694{{u_i}}^{2}{q}^{3}\right){x_3}{\Delta_1}
$$ $$
+\left(
-143992{q}^{3}
-46818{u_i}{q}^{2}
+ 4522{{u_i}}^{2}q\right){x_2}{\Delta_2}
$$ $$
+\left(
-491334{q}^{4}
-164832{u_i }{q}^{3}
+15822{{u_i}}^ {2}{q}^{2}\right){x_3}{\Delta_2}
$$ $$
+\left(
35042{q}^{2}
+11348{u_i}q
-1098{{u_i}}^{2}\right){x_2}{\Delta_3}
$$ $$
+\left(
112272{q}^{3 }
+ 37824{u_i}{q}^{2}
-3633{{u_i}}^{2}q\right){x_3} {\Delta_3}.
$$

>From this we compute the special coordinates
$$
\eta_0 = -\frac{1}{2q^3}, \quad \quad
\eta_i = 49346q +16192{u_i} -1562{\frac {{{u_i}}^{2}}{q}},
$$ $$
\eta_{00} = -\frac{1}{q^4}, \quad \quad
\eta_{i0} = \eta_{0i} = 
-{\frac {-2536 q^2-688\,{u_i} q+69\,{{u_i}}^{2}}
{968\,{q}^{6}}},
$$ $$
\eta_{ii} =
\frac {-3412 q^2-260\,{u_i}q+41\,{{u_i}}^{2}}{484\,q^{6}}
$$
$$
\eta_{ij} = 
{\frac {404}{121\,q^4}}
-{\frac {34}{121}}\,{\frac {u_{i} + u_j}{{q}^{5}}}
+{\frac {13}{968}} \,{\frac {u_{{j}}u_{{i}}}{{q}^{6}}}.
$$

\newpage

\centerline{\bf \S 3. Del Pezzo surfaces}

\medskip

{\bf 3.1. Notation.}
Let $V$ be a surface. For surfaces,
we will denote the dual class of a point by $\Delta_2$, and
the dual coordinate by $z$.  Assuming $h^{2,0}(V)=0,$
we may and will identify the group
$H_2(V,\bold{Z})/(tors)$ with
$H^{1,1}(V,\bold{Z})$.

\smallskip

We start again with calculating the structure constants
of the quantum multiplication of classes in $H^{p,p}$ restricted to the first
infinitesimal neighborhood of $H^{1,1}$. Essentially,
this means that we must calculate the third
derivatives of the potential modulo $z^2$. 
The canonical coordinates of this ``very small quantum cohomology''
are the eigenvalues of the operator $k\circ ,\, k=c_1(V),$
calculated at the points of $H^{1,1}$ so that
only triple correlators are involved in their calculation.

\medskip

{\bf 3.2. Quantum multiplication table.} Since $\Delta_0$
is the identity with respect to $\circ$, we need only the following
formulas, which follow from (0.1). Let $D_1,\,D_2\in H^2(V),$ then
$$
D_1\circ D_2= (D_1,D_2)\,\Delta_2 \,+
$$
$$
\sum_{\beta \ne 0}\langle e^{z\Delta_2}\rangle_{\beta}
(D_1,\beta )\,(D_2,\beta )\,q^{\beta}\beta \,+
\left(\sum_{\beta \ne 0}
\langle \Delta_2e^{z\Delta_2}\rangle_{\beta}(D_1,\beta )\,
(D_2,\beta )q^{\beta}\right) \Delta_0 .
\eqno(3.1)
$$

Similarly, for $D\in H^2(V)$,
$$
D\circ\Delta_2 =
\sum_{\beta \ne 0}\langle\Delta_2 e^{z\Delta_2}\rangle_{\beta}
(D,\beta )\,q^{\beta}\beta \,+
\left(\sum_{\beta \ne 0}
\langle \Delta_2^2e^{z\Delta_2}\rangle_{\beta}(D,\beta )\,
q^{\beta}\right)\, \Delta_0 .
\eqno(3.2)
$$

Finally,
$$
\Delta_2\circ\Delta_2 =
\sum_{\beta \ne 0}\langle\Delta_2^2e^{z\Delta_2}\rangle_{\beta}
\,q^{\beta}\beta \,+
\left(\sum_{\beta \ne 0}
\langle \Delta_2^3e^{z\Delta_2}\rangle_{\beta}\,
q^{\beta}\right)\, \Delta_0 .
\eqno(3.3)
$$
If we want only terms modulo $z^2$, (1.5) can be used
in order to restrict the summation over $\beta$:
$$
D_1\circ D_2= (D_1,D_2)\,\Delta_2 \,+
$$
$$
\sum_{\beta :\,k(\beta )=1}\langle \emptyset\rangle_{\beta}
(D_1,\beta )\,(D_2,\beta )\,q^{\beta}\beta +
\left(\sum_{\beta :\,k(\beta )=2}
\langle \Delta_2\rangle_{\beta}(D_1,\beta )\,
(D_2,\beta )q^{\beta}\right)\, \Delta_0 \,+
$$
$$
z\sum_{\beta :\,k(\beta )=2}\langle \Delta_2\rangle_{\beta}
(D_1,\beta )\,(D_2,\beta )\,q^{\beta}\beta \,+\,
z\,\left(\sum_{\beta :\,k(\beta )=3}
\langle \Delta_2^2\rangle_{\beta}(D_1,\beta )\,
(D_2,\beta )q^{\beta}\right)\, \Delta_0 \,+\, O(z^2),
\eqno(3.4)
$$
$$ 
D\circ \Delta_2 =
\sum_{\beta :\,k(\beta )=2}\langle \Delta_2\rangle_{\beta}
(D,\beta )\,q^{\beta}\beta +
\left(\sum_{\beta :\,k(\beta )=3}
\langle \Delta_2^2\rangle_{\beta}(D,\beta )\,
q^{\beta}\right)\, \Delta_0 \,+
$$
$$
z\sum_{\beta :\,k(\beta )=3}\langle \Delta_2^2\rangle_{\beta}
(D,\beta )\,q^{\beta}\beta +
z\left(\sum_{\beta :\,k(\beta )=4}
\langle \Delta_2^3\rangle_{\beta}(D,\beta )\,
q^{\beta}\right)\, \Delta_0 \,+\, O(z^2).
\eqno(3.5)
$$

In the case when one of the divisor classes is $k$, (3.4) and (3.5)
further simplify:
$$
k\circ D= k(D)\,\Delta_2 \,+
$$
$$
\sum_{\beta :\,k(\beta )=1}\langle \emptyset\rangle_{\beta}
(D,\beta )\,q^{\beta}\beta +
2\left(\sum_{\beta :\,k(\beta )=2}
\langle \Delta_2\rangle_{\beta}(D,\beta )\,
q^{\beta}\right)\, \Delta_0 \,+
$$
$$
2z\sum_{\beta :\,k(\beta )=2}\langle \Delta_2\rangle_{\beta}
(D,\beta )\,q^{\beta}\beta \,+\,
3z\,\left(\sum_{\beta :\,k(\beta )=3}
\langle \Delta_2^2\rangle_{\beta}(D,\beta )\,
q^{\beta}\right)\, \Delta_0 \,+\, O(z^2).
\eqno(3.6)
$$

In particular,
$$
k\circ k= (k,k)\,\Delta_2 \,+
\sum_{\beta :\,k(\beta )=1}\langle \emptyset\rangle_{\beta}
\,q^{\beta}\beta +
4\left(\sum_{\beta :\,k(\beta )=2}
\langle \Delta_2\rangle_{\beta}\,
q^{\beta}\right)\, \Delta_0 \,+
$$
$$
4z\sum_{\beta :\,k(\beta )=2}\langle \Delta_2\rangle_{\beta}
\,q^{\beta}\beta \,+\,
9z\,\left(\sum_{\beta :\,k(\beta )=3}
\langle \Delta_2^2\rangle_{\beta}\,
q^{\beta}\right)\, \Delta_0 \,+\, O(z^2).
\eqno(3.7)
$$

Furthermore,
$$ 
k\circ \Delta_2 =
2\sum_{\beta :\,k(\beta )=2}\langle \Delta_2\rangle_{\beta}
\,q^{\beta}\beta +
3\left(\sum_{\beta :\,k(\beta )=3}
\langle \Delta_2^2\rangle_{\beta}\,
q^{\beta}\right)\, \Delta_0 \,+
$$
$$
3z\sum_{\beta :\,k(\beta )=3}\langle \Delta_2^2\rangle_{\beta}
\,q^{\beta}\beta +
4z\left(\sum_{\beta :\,k(\beta )=4}
\langle \Delta_2^3\rangle_{\beta}\,
q^{\beta}\right)\, \Delta_0 \,+\, O(z^2).
\eqno(3.8)
$$
\medskip

{\bf 3.3. Del Pezzo surfaces.} Let $V=V_r$ be obtained from $\bold{P}^2$ by blowing up $r$ points
in general position. We get del Pezzo surfaces (characterized by the
ampleness of $k$) for $r\le 8.$
Sums in the right hand sides of
(3.4)--(3.8) for larger ranks $r+1$ of $H^{1,1}(V)$ are long
and for $r\ge 9$ are infinite. However, for del Pezzo surfaces
(and to a certain degree, for more general rational surfaces)
the situation becomes simpler because a Weyl group
acts upon $\beta \in H^{1,1}(V)$ leaving $k(\beta ),\,(\beta ,\beta )$
and the correlators 
$\langle\Delta_2^{k(\beta )-1}\rangle_{\beta} $ invariant. For $k(\beta )\le 4$ the number of orbits is quite small, and the resulting expressions
become manageable. We start with some classical results
on the structure of $H^2(V,\bold{Z})$ (see e.~g. [Ma3]).

\smallskip

Fix a representation of $V_r$ as blow up of $\bold{P}^2.$
Then $k$ and the classes  of the resulting
exceptional curves $l_1,\dots ,l_r$ generate a sublattice of index
$3$ in the N\'eron--Severi group coinciding with $H^2(V,\bold{Z}).$
The whole lattice $N_r$ is generated by the orthonormal basis
$(l_0;\,l_1,\dots ,l_r)$ where $l_0$ is the class of a line lifted
to $V.$ We identify $\beta =al_0-\sum_ib_il_i$ with
the vector $(a;\,b_1,\dots ,b_r)$ so that $a=a(\beta )=(\beta ,l_0)$
and $b_i=b_i(\beta )=(\beta, l_i).$ We have $k=k_r=(3;1,\dots ,1)$
in this notation.

\smallskip

Denote by $R_r$ the set $\{\beta\,|\,k(\beta )=0,\,(\beta ,\beta )=-2\}.$
It is a system of roots, lying entirely
in the orthogonal complement $O_r$ to $k_r$. Reflections with respect to
the roots generate the Weyl group $W_r$ fixing $k_r$. In our representation,
this group for $r\ge 3$ is generated by all permutations of $l_r$'s and
by the Cremona transformation
$$
(a;\,b_1,\dots ,b_{r})\mapsto (a+\delta ;\, b_1+\delta, b_2+\delta, b_3+\delta,
b_4, \dots ,b_{r}),\quad \delta := a-b_1-b_2-b_3.
\eqno(3.9)
$$ 
\smallskip

The function
$\beta\mapsto \langle\Delta_2^{k(\beta)-1}\rangle_{\beta}$
is constant on $W_r$--orbits. This is well known. One argument
is given in [GP], 5.1. It involves geometric Cremona
transformations (blowing up and down), which generate $W_r$, together
with permutations of $(l_1,\dots ,l_r)$. Alternatively,
one can invoke
Hirschowitz's theorem ([Hirsch]) which shows that
the action of $W_r$ on $N_r$ is induced by the action
of $W_r$ on a generic del Pezzo surface and its definition
field. 

\smallskip

In fact, a related reasoning readily extends 
to a more general situation. Consider a family of
algebraic manifolds over a base whose fundamental
group acts on the fiber cohomology in a nontrivial way.
A typical situation is the existence of several
Lefschetz type transforms associated with various
cusps of the moduli space. If they are reflections,
a version of the following proposition will hold.

\medskip

\proclaim{\quad 3.3.1. Proposition} Let $\rho$ be a root in
$R_r$, $s_{\rho}: \gamma\mapsto \gamma +(\gamma ,\rho )\,\rho$ the respective reflection. It induces the involution of the algebraic torus $T_{V_r}$ with the character group
$H^2(V_r,\bold{Z})$. Let $T_{\rho}$ be the codimension one
subtorus consisting of fixed points of this involution:
$q^{s_{\rho}(\beta )} = q^{\beta}$ on $T_{\rho}$ for all $\beta$.
Then we have, after restriction on $T_{\rho}$:

\smallskip

a) For any $D\in H^2(V)$ orthogonal to $\rho$
$$
D\circ\rho =\mu_{D,\rho}\,\rho ,\ \mu_{D,\rho}=
-\frac{1}{2}\sum_{\beta} \langle e^{z\Delta_2}\rangle_{\beta}\,(\rho ,\beta )^2
(D,\beta )q^{\beta} .
\eqno(3.11)
$$
Similarly,
$$
\Delta_2\circ\rho =\nu_{\rho}\,\rho ,\  \nu_{\rho}=
-\frac{1}{2}\sum_{\beta} \langle e^{z\Delta_2}\rangle_{\beta}\,(\rho ,\beta )^2
\,q^{\beta} .
\eqno(3.12)
$$

b) If $D_1, D_2$ are orthogonal to $\rho$, then $D_1\circ D_2$
is orthogonal to $\rho$ on $T_{\rho}$; the same is true for $\rho\circ\rho$.

\smallskip
\endproclaim

\smallskip

{\bf Proof.} We have the following general identity. Let $S$
be a $W_r$--orbit, $D_1, \dots ,D_n$ orthogonal to $\rho$,
and $m$ an integer. Then we have on $T_{\rho}$
$$
\sum_{\beta\in S} (D_1,\beta )\dots (D_n,\beta )\,(\rho ,\beta )^{2m+1}
q^{\beta} =0\,.
\eqno(3.13)
$$
To see this, replace $s$ by $s_{\rho}(\beta )$ in the summands
of (3.13). The whole sum will remain the same,
whereas each summand will change sign, because of the
factor $(\rho ,\beta )^{2m+1}$.

\smallskip

Consider now (3.1) for $D_1=D$, $D_2=\rho$. In view of
(3.13), the coefficient of $\Delta_0$ vanishes on $T_{\rho}.$
Take the scalar product with any $D^{\prime}$ orthogonal to
$\rho$. It will vanish, again in view of (3.13). Hence
our vector is proportional to $\rho$. To calculate
the coefficient $\mu_{\rho}$, it suffices to take
the scalar product with $\rho$. In order to prove (3.12),
we treat similarly (3.2) with $D=\rho$.  

\smallskip

To establish b), look at (3.1) with both $D_i$
orthogonal to $\rho$, resp. with $D_1=D_2=\rho$.
(Notice that in these cases the coefficient of $\Delta_0$
does not vanish, but of course, $\Delta_0$ is still
orthogonal to $\rho$ in the Poincar\'e metric).
This ends the proof.

\medskip

The intersection of all subtori
$T_{\rho}$ is the one--dimensional subgroup
$q^{\beta}=e^{tk_r(\beta )}$, $t$ arbitrary. On this subgroup,
and in particular at $t=0$ (which we write as $q=1$),
the statements of the Proposition  3.3.1 are valid
simultaneously for all roots $\rho$. Hence we get the following corollary.

\medskip

\proclaim{\quad 3.3.2. Corollary} On the subgroup defined above,
$k_r\circ$, $\Delta_2\circ$ act on the orthogonal complement $O_r$ to $k_r$
as multiplication by the scalars (independent on the choice
of $\rho\in R_r$) denoted respectively
$$
\mu = -\frac{1}{2}\sum_{\beta}\langle e^{z\Delta_2}\rangle_{\beta}\,
(\rho ,\beta)^2\,k_r(\beta )\,e^{tk_r(\beta )},
\eqno(3.14)
$$
$$
\nu = -\frac{1}{2}\sum_{\beta}\langle \Delta_2 e^{z\Delta_2}\rangle_{\beta}\,
(\rho ,\beta)^2\,e^{tk_r(\beta )},
\eqno(3.15)
$$
\endproclaim

\smallskip

In fact, $O_r$ is spanned by all roots in $R_r$.

\medskip

{\bf 3.4. Structure of the sets $k(\beta )=a.$} For each
$3\le r\le 8$ we denote by $I_r, F_r, G_r, H_r$ the following subsets of $N_r$,
each constituting one orbit of $W_r$. The statements accompanying
their definitions are proved, for example, in [Ma3], Ch.~4.

\smallskip

$$
I_r:=\{\beta\,|\,k(\beta )=(\beta,\beta )=-1\}=W_r(0;-1;0,\dots ,0).
\eqno(3.16)
$$
Elements of $I_r$ are exactly classes of exceptional curves on a general
del Pezzo surface of degree $9-r.$ (This is also true for $r=9$
but false for larger values of $r$: cf. [Hirsch], 3.4).
$$
F_r:=W_r(l_0-l_1)=W_r(1;1,0,\dots ,0).
\eqno(3.17)
$$
$$
G_r:=W_rl_0=W_r(1;0,0,\dots ,0).
\eqno(3.18)
$$
$$
H_r:=W_r(l_1+l_2)=W_r(0;-1,-1,0, \dots ,0).
\eqno(3.19)
$$
Elements of $H_r$ for $r\ge 4$ are exactly classes
of cycles $\lambda +\mu$ where $\lambda$ and $\mu$ are
two disjoint exceptional curves. For $r=3$, however, there are two
$W_3$--orbits of such cycles: blowing down  two
disjoint exceptional curves can produce either $\bold{P}^2$ blown up
at a single point or $\bold{P}^1\times\bold{P}^1.$ The orbit $H_3$
corresponds to the first case.

\medskip

\proclaim{\quad 3.4.1. Proposition} The total support and
nonzero values of the function 
$\beta\mapsto \langle\Delta_2^{k(\beta)-1}\rangle_{\beta}$ for $k(\beta )\le 3$ are given in the following list:
\smallskip

a) $\langle\emptyset\rangle_{\beta}=1$ on $I_r$ for all
$3\le r\le 8.$

In $N_8$ in addition $\langle\emptyset\rangle_{k_8}=12$ (recall that $k_r$ is
$W_r$--invariant for any $r$).

\smallskip

b) $\langle\Delta_2\rangle_{\beta}=1$ on $F_r$ for all 
$3\le r\le 8$. 

In $N_7$ in addition $\langle\Delta_2\rangle_{k_7}=12.$

In $N_8$ in addition $\langle\Delta_2\rangle_{\beta}=12$ 
on $k_8 + I_8$ and $\langle\Delta_2\rangle_{2k_8}=90.$

\smallskip

c) $\langle\Delta_2^2\rangle_{\beta}=1$ on $G_r$ for all 
$3\le r\le 8$. 

In $N_6$ in addition $\langle\Delta_2^2\rangle_{k_6}=12.$

In $N_7$ in addition $\langle\Delta_2^2\rangle_{\beta}=12$
on $k_7+I_7.$

In $N_8$ there are four additional orbits with the
following values: $12$ on $k_8+H_8$;
$96$ on $k_8+F_8$; $576$ on $2k_8+I_8$; $2880$ on
$3k_8$.  
\endproclaim

\smallskip

{\bf Proof.} 
The argument depends on compiling
and studying a table, fortunately rather short one.

\smallskip

It lists all non--zero 9--uples $\beta = (a;b_1,\dots , b_8)$ 
in $\bold{Z}^9$ with the
following properties:
$$
1\le k(\beta )=3a-b_1-\dots -b_9\le 3,\ b_1\ge b_2\ge \dots \ge b_8\ge 0,\
a\ge b_1+b_2 +b_3. 
\eqno(3.20)
$$
In addition, it includes $(0;0,-1,0,\dots ,0)$ and contains
10 entries altogether. In the notation omitting final zeros
and rendering multiplicities of $b_i$'s by superscripts (cf. [GP], p. 87), they are
$$
(0;-1),\, (1;),\,(1;1),\, (3;1^6),\,(3;1^7),\, (3;1^8),\,(4;2,1^7),\,
(6,2^7,1),\, (6;2^8),\, (9;3^8)\, .
\eqno(3.21)
$$

According
to [GP], sections 3 and 5, any $\beta\in N_8$ with  $\langle\Delta_2^{k(\beta)-1}\rangle_{\beta} \ne 0$ and $k(\beta )\le 3$
is $W_8$--equivalent to one of the entries of this list.
(Actually, to a single entry, because they all are pairwise distinguished by the
values of pairs $k(\beta ), (\beta ,\beta )$). 
Moreover, if an entry ends with $\ge s$ zeroes,
then for $s\ge 3$ its $W_{8-s}$--orbit contains
all $\beta\in N_{8-s}$ with non--zero correlators
and respective value of $k(\beta ).$ 
\smallskip

In fact, for any such $\beta$
we can
first permute $b_i$'s to make them decreasing. If we get
then  $\delta :=a -b_1-b_2-b_3<0$, we can decrease $a$ by applying
the respective Cremona transformation  (3.9)
or, which is the same,
by reflecting $\beta$ with respect of one of the simple roots
in $N_{8-s}.$ If after several steps of this kind we get
a vector with $a>0$ and some negative $b_i$, or with $a=0$
but not with exactly one $b_i=-1$ and zeros on the remaining
places, then $\langle\Delta_2^{k(\beta)-1}\rangle_{\beta} = 0$ 
according to [GP]. Hence we have to consider only
the 10 classes listed above. 

\smallskip

The values of $\langle\Delta_2^{k(\beta)-1}\rangle_{\beta}$ for these classes
can be read off from the tables on p. 88 of [GP]. One exception
is $(9;3,\dots ,3)$ which must be calculated using recursion;
its value was communicated to us by L.~G\"ottsche.

\smallskip

Finally, we can directly identify the relevant
$W_r$--orbits with the sets listed in the Proposition.
Notice however that adding zeros changes the
invariant description of the vector. For example,
$(3;1^6)$ produces $k_6$ in $N_6$, but $k_7+l_7$
in $N_7$, and for the respective orbits we have
$W_6k_6=\{k_6\},\,W_7(k_7+l_7)=k_7+I_7$.

\smallskip

The complete list of $F_8$ (calculated for other purposes)
in this format can be found in [MaTschi] on p\. 329.
Similarly, the complete list of $G_8$ can be extracted
from the table on p\. 330 in the following way:
replace the misprinted entry $(10;5^2,4^3,3^2,2)$
by the correct one $(10;5^2,3^5,2)$ and delete
the following six entries belonging to another
orbit (with vanishing Gromov--Witten invariant):
$(5;3^2,1^6),\,(7;3^5,1^3),\,(9;5,3^6,1),\,
(11;5^3,3^5),\,(13;5^6,3^2),\,(15;7,5^7).$

\medskip

{\bf 3.5. The point $q=1$, $z=0$.} 
We will state some simple lemmas which use the $W_r$-invariance
and will simplify the calculation of the quantum cup product; at
the end of this section we will combine the obtained information to 
describe the quantum cohomology algebra at $q=1$, $z=0$ in the
case $r \ge 5$.

\smallskip

First, note the following

\proclaim{\quad 3.5.1 Fact} If $r \ge 5$ and $\rho$ is a root, then
the orthogonal subspace $\rho^\bot$ inside the root space $O_r$ is
generated by the roots orthogonal to $\rho$.
\endproclaim

\smallskip

Indeed, if $r=5$ and $\rho = (0; -1, 1,0,0,0)$, then $(1; 1,1,1,0,0)$,
$(1;0,0,1,1,1)$,
$(0; 0,0,-1, 1, 0)$ and $(0; 0,0,0,-1,1)$ are linearly independent
roots and all orthogonal to $\rho$. This argument can easily be extended
to higher $r$.

\smallskip

This fact will actually make a difference in several computations, which is
why the Proposition 3.5.5 would be false for $r \le 4$.

\medskip

\proclaim{\quad 3.5.2. Lemma} Let $S$ be a $W_r$ orbit
in $N_r$, $k_r(S):=k_r(\beta )$ for an arbitrary $\beta \in S$, and
let $\rho \bot \rho'$ be two orthogonal roots.
Then we have:

\smallskip

a) $\sum_{\beta\in S}\beta = k_r(S)\,\dfrac{|S|}{9-r}\,k_r.$

\smallskip

b) $\sum_{\beta \in S} (\rho, \beta) \beta
		= -\frac{1}{2} \sum_{\beta \in S} (\rho, \beta)^2 \rho$.

\smallskip

c) $\sum_{\beta \in S} (\rho, \beta)^{m} (\rho', \beta) = 
 \sum_{\beta \in S} (\rho, \beta)^{2m+1} = 0$.

\smallskip

d) If $r \ge 5$, then
	$\sum_{\beta \in S} (\rho, \beta)^2 \beta
		= \dfrac{k_r(S)}{9-r} \sum_{\beta \in S} (\rho, \beta)^2\,k_r$.

\smallskip

e) $ \sum_{\beta \in S} (\rho, \beta) (\rho',\beta) \beta = 0$.
\endproclaim

\smallskip

{\bf Proof.}
a) Since  $\sum_{\beta\in S}\beta$ is a $W_r$--invariant
element of $N_r$, it must be proportional to $k_r$.
To calculate the proportionality coefficient, it remains to take
the intersection index with $k_r$.

\smallskip

b) This follows from symmetrizing over $\beta$ and
$s_\rho(\beta) = \beta + (\rho,\beta) \rho$:
$$
\sum_{\beta \in S} (\rho, \beta) \beta =
\frac{1}{2} \sum_{\beta \in S}
	\Bigl( (\rho, \beta) \beta + (\rho, s_\rho(\beta)) s_\rho(\beta)
								\Bigr)
	=
\frac{1}{2} \sum_{\beta \in S}
	\Bigl( (\rho, \beta) \beta - (\rho, \beta) (\beta + (\rho,\beta)\rho)
								\Bigr).
$$

\smallskip

c) These are special cases of (3.13).

\smallskip

d) Since we assume $r \ge 5$, this formula can be checked by taking
the intersection index of the right- and left-hand side with $k_r$, with
$\rho$ and with an arbitrary root $\rho' \bot \rho$, respectively.
For $\rho, \rho'$ this follows from c), and for $k_r$ it is obvious.

\smallskip

e) As in b), we symmetrize over $\beta$ and $s_\rho(\beta)$, from which
we see that the sum is equal to
$-\frac{1}{2} \sum_{\beta \in S} (\rho, \beta)^2 (\rho',\beta) \rho$;
this is zero according to c).

\medskip

{\bf 3.5.3. Cardinalities.} The following table lists
the cardinalities of the sets introduced above:

\medskip

\settabs 7\columns
\+ $r$ & 3 & 4& 5& 6& 7& 8 \cr

\+ $|W_r|$ &$2^23$ & $2^33\cdot 5$& $2^73\cdot 5$& $2^73^45$& 
$2^{10}3^45\cdot 7$& $2^{14}3^55^27$ \cr 

\+ $|R_r|$ &$2^3$ & $2^2\cdot 5$& $2^3\cdot 5$& $2^33^2$& 
$2\cdot 3^27$& $2^{4}3\cdot 5$ \cr 

\+ $|I_r|$ &$2\cdot 3$&$2\cdot 5$ & $2^4$& $3^3$& $2^37$& 
$2^43\cdot 5$ \cr

\+ $|F_r|$ & 3&$ 5$ & $2\cdot 5$& $3^3$& $2\cdot 3^27$& 
$2^43^35$ \cr

\+ $|G_r|$ & 2&$ 5$ & $2^4$& $2^33^2$& $2^63^2$& 
$2^73^35$ \cr

\+ $|H_r|$&$2\cdot 3 $&$2\cdot 3\cdot 5 $&$2^4\cdot 5$ &
$2^33^3$ &$2^23^37$& $2^63\cdot 5\cdot 7$\cr

\medskip

The orders of $W_r,\,R_r$ and $I_r$ are well known. The remaining
entries can be calculated directly for $r=3$ and
then inductively in the following way.
There are exactly $r$ exceptional classes in $N_r$
having zero intersection index with $(1;0^r)$, and
they are pairwise orthogonal. Hence $|G_r|$ contains as many
elements as there are maximal contractible $r$--tuples of 
exceptional classes in $N_r$, that is $|I_r|\,|G_{r-1}|/r.$
Similarly, $|F_r|$ counts non--maximal contractible
($r-1$)--tuples so that $|F_r|=|I_r|\,|F_{r-1}|/2(r-1).$
Finally, $|H_r|$ counts unordered contractible pairs, 
contained in maximal contractible $r$--tuples so that
$|H_r|=|I_r|\,|I_{r-1}|/2$.

\medskip

\proclaim{\quad 3.5.4. Proposition} At $q=1,\,z=0,$ the
operator $k_r\circ$ admits two complementary orthogonal invariant
subspaces: one spanned by $\Delta_0,k_r,\Delta_2,$
the other being $O_r$.

\smallskip

On the first subspace, the characteristic polynomial 
$\roman{det}\,(\lambda\,\roman{id} -k_r\circ )$ is
$$
R(\lambda) = \lambda^3 -B_r\lambda^2 -D_r\lambda -36\,(9-r)\,C_r
\eqno(3.22)
$$
where
$$
B_r=\frac{|I_r|}{9-r}\,,\ C_r=\frac{|G_r|}{12} +\delta_{r,6}
+56\,\delta_{r,7}+35760\,\delta_{r,8}, 
$$
$$
D_r=8\,(|F_r|+12\,\delta_{r,7}+12\,\delta_{r,8}|I_8|
+90\,\delta_{r,8})
\eqno(3.23)
$$
($\delta_{r,a}$ being the Kronecker symbol).
\smallskip

On the second subspace, $k_r\circ$ is multiplication by a constant
$\mu_r$ if $r \ge 4$. This constant is given by the following table:

\medskip

\settabs 6\columns
\+ $r$ & $4$& $5$& $6$& $7$& $8$ \cr

\+ $\mu_r$ & $-3$& $-4$& $-6$& $-12$& $-60$ \cr

\endproclaim 

\smallskip

{\bf Proof.} Formulas (3.7) and (3.8), combined
with Proposition 3.4.1 and Lemma 3.5.1, lead to the following
action of $k_r\circ$ at $q=1,\,z=0$ on the first subspace:
$k_r\circ\Delta_0=k_r$ and
$$
k_r\circ k_r=4\,\bigl(|F_r|+12\,\delta_{r,7}+12\,\delta_{r,8}|I_8|
+90\,\delta_{r,8}\bigr)\Delta_0
+\left(\frac{|I_r|}{9-r}+12\,\delta_{r,8}\right)k_r + (9-r)\,\Delta_2,
					\eqno(3.24)
$$
$$
k_r\circ\Delta_2=2\left(\frac{2\,|F_r|}{9-r} +12\,\delta_{r,7} +24\,\delta_{r,8}|I_8|
+180\,\delta_{r,8}\right)\,k_r +
					\eqno(3.25)
$$
$$
3\,\bigl(\,|G_r|+ 12\,\delta_{r,6}+12\,\delta_{r,7}|I_7|
+\delta_{r,8}(12\,|H_8|+96\,|F_8|+576\,|I_8|+2880\,)\,\bigr)\Delta_0 .
$$
The direct calculation of the characteristic polynomial then gives
(3.22) and (3.23).
\smallskip

To treat the orthogonal complement $O_r$ of $k_r$ in $H^2(V)$, notice
first of all that it is spanned
by roots $\rho \in R_r$. If $r\ge4$, they form a single
$W_r$--orbit. We calculate $k_r\circ\rho$ at $q=1,z=0$ using (3.6).
The coefficient of $\Delta_0$ vanishes according to 3.5.2 c),
while the first sum in (3.6) is equal to
$-\frac{1}{2}\suma (\rho,\beta)^2 \rho =: \mu_r$ by
lemma 3.5.2 b). In view of $W_r$--invariance, the coefficient cannot depend on
the root $\rho$ if $r \ge 4$.

\smallskip

The value of $\mu_r$ can be computed explicitly by using a complete list
of the $W_r$--orbit $I_r$ as given in [Ma3]; see the relevant remark in
the proof of the Proposition 3.5.5 for a more elegant way to compute it.

\medskip

\proclaim{\quad 3.5.5. Proposition} If $r\ge 5$, then the quantum cohomology
algebra at the point $q=1, z=0$ is isomorphic to
$ \C \oplus \C [X_1, \dots, X_r] /(X_i^2=X_j^2, X_iX_j = 0)$. To
establish this isomorphism, we can send $X_1, \dots, X_r$ to an arbitrary
orthonormal basis of $O_r$, in which case $X_1^2 = \dots = X_r^2$
gets sent to
$$
\Gamma := \Delta_2 
 	  - \frac{1}{2(9-r)} \suma (\rho, \beta)^2 k_r
	  - \frac{1}{2} \sumb (\rho, \beta)^2 \Delta_0 .
								\eqno(3.26)
$$
\endproclaim

{\bf Proof.} Let $\rho$ be a root. Using lemma 3.5.2 d) we compute
$$
\rho \circ \rho
	= - 2 \Delta_2 + \suma (\rho, \beta)^2 \beta
		+ \sumb (\rho, \beta)^2 \Delta_0 =
$$
$$
 	 -2 \Delta_2 + \frac{1}{9-r} \suma (\rho, \beta)^2 k_r
		+ \sumb (\rho, \beta)^2 \Delta_0 = -2 \Gamma.
$$

Now let $\rho'$ be a root orthogonal to $\rho$. By 3.5.2 c) and e) we
see that
$$
\rho \circ \rho' =
	  \suma (\rho, \beta) (\rho',\beta) \beta
	+ \sumb (\rho, \beta) (\rho',\beta) \Delta_0
	= 0.
$$

Since $\rho$ and the roots $\rho'$ orthogonal to $\rho$ span $O_r$,
we can conclude that for any $D'$ which lies in this space, we have
$ \rho \circ D'  = (\rho, D')\,\Gamma$. Since the roots span $O_r$,
the formula $ D \circ D' = (D, D')\, \Gamma$ holds for
all $D, D' \in O_r$.

\smallskip

Note that $\Gamma \circ \rho = \rho' \circ \rho' \circ \rho = 0$ and
hence $\Gamma \circ \Gamma = 0 = \rho^{\circ 3}$.

\smallskip

By associativity, $\Gamma $ must be an eigenvector of $\adkr$ to the
eigenvalue $\mu_r$. So $\mu_r$ must be a root of the characteristic
polynomial $R(\lambda)$ defined in (3.22). One can check that it is
actually a double root. In fact, if it was a simple
root, then the eigenspace of $\adkr$ to the eigenvalue $\mu_r$ would
consist only of nilpotent elements, which is impossible since it is
a direct summand of an algebra with identity. (This fact gives the most
direct way to determine $\mu_r$: as the double root of $R(\lambda)$.)
Since $R(0)<0$, it can't be a triple root.

\smallskip

Thus we have an $(r+2)$--dimensional eigenspace of $\adkr$ and a splitting
of $H^*(V)$ into a one--dimensional and an $(r+2)$--dimensional subalgebra. 
The latter one has a basis consisting of its identity, $\Gamma$ and an
arbitrary orhonormal basis of $k_r^\bot$; in this basis the structure
is evidently as described in the proposition. 

\medskip

{\bf 3.6. Semisimplicity of the quantum cohomology of del Pezzo surfaces.}
We have seen that if $r \ge 5$, the quantum cohomology is far from being
semisimple at the symmetric point $q=1$; also proposition 3.5.4 shows that
this cannot be a tame semisimple point if $r \ge 3$. However, we will prove:

\smallskip

\proclaim{\quad 3.6.1. Theorem} The Frobenius manifold associated to the
quantum cohomology of a del Pezzo surface is generically semisimple on
$H^2$.
\endproclaim

\smallskip

Our proof is based on extending the Frobenius manifold associated to $V_r$
to a boundary; on the boundary, it ceases to be a Frobenius manifold, yet
the boundary itself will be the Frobenius manifold associated to $V_{r-1}$.

\smallskip

First, we introduce appropriate coordinates:
We write the generic character $q^\beta$ on $H^2$ as
$$
q^\beta = e^{\sum_{i=0}^{r} x_i (l_i, \beta)},
$$
where the $l_i$ are as defined in the beginning of section 3.3. We define
$q_i := e^{x_i}$ and thus can write
$$
q^\beta = \prod_{i=0}^r q_i^{(l_i,\beta)}.
$$

As mentioned in the beginning of section 3.3, all the sums appearing in
equations (3.4) to (3.6) are finite, and each summand is a monomial in the
variables $q_0, \dots, q_r$ and $q_1^{-1}, \dots, q_r^{-1}$. Hence, the
spectral cover map restricted to $T_V$ is the morphism of affine schemes  
$$
\Spec\,H^*(V_r,\Q)[q_0, q_1, \dots, q_r, q_1^{-1}, \dots q_r^{-1}] \to 
\Spec\,\Q[q_0, q_1, \dots, q_r, q_1^{-1}, \dots q_r^{-1}] 	\eqno(3.27)
$$
where the left hand side is the small quantum cohomology algebra.
\smallskip

\proclaim{\quad 3.6.2. Lemma} The map (3.27) can be extended in a flat way
to the boundary $q_r = 0$.
Moreover, the fiber over $q_r=0$ is isomorphic to a disjoint sum of the
identity map and the respective map (3.27) for the case $r-1$.
\endproclaim

\smallskip

{\bf Proof.} As we already mentioned,
from the properties (P1) and (P2) of Theorem 4.1 in [GP] it follows that
the only $\beta \in H_2 (V,\Z)$ with non--vanishing Gromov--Witten invariant
$\langle \Delta_2^{k_r(\beta)-1} \rangle$ and $(\beta, l_r) < 0$ is $l_r$
itself, for which we have $\langle \emptyset \rangle_{l_r} = 1$.
This will be the key fact used in our proof.

\smallskip

Let $F$ be the free module over
$\Q[q_0, \dots, q_r, q_1, \dots, q_{r-1}^{-1}]$ with basis
$\Delta_0, \Delta_2, \allowbreak{} l_0, \dots, \allowbreak{} l_{r-1},
q_r l_r$. The module $F[q_r^{-1}]$ is simply the small
quantum cohomology algebra; we will show that $F$ is a subalgebra
by simply examining all relevant products.

\smallskip

Let $D$ be an arbitrary element in $H^2(V_r)$ orthogonal to $l_r$, i.~e.
$D \in H^2(V_{r-1})$, and let $T_1$ and $T_2$ be two arbitrary elements in
$H^*(V_r)$ orthogonal to $l_r$, so $T_1, T_2 \in H^*(V_{r-1})$. We will
denote their product by $\circ_{V_r}$ or $\circ_{V_{r-1}}$ indicating in
which quantum cohomology this product is to be taken.

\smallskip

We claim that:
$$
q_r l_r \circ q_r l_r = q_r l_r + O(q_r^2)
$$
$$
q_r l_r \circ D = O(q_r^2)
$$
$$
q_r l_r \circ \Delta_2 = O(q_r^2)
$$
$$
T_1 \circ_{V_r} T_2 = T_1 \circ_{V_{r-1}} T_2 + O(q_r)		\eqno(3.27)
$$

All these formulas follow from (3.4)--(3.6) and the
fact mentioned in the beginning of our proof. To give an example of the
computations, we treat the case of $q_r l_r \circ D$:
$$
l_r \circ D = 
\sum_{\beta:\,k_r(\beta)=1} \langle \emptyset \rangle_\beta
	q_0^{(l_0,\beta)} \cdot \dots\cdot 
	q_r^{(l_r,\beta)} (D,\beta) (l_r, \beta) \beta
$$
$$
+ \sum_{\beta:\,k_r(\beta)=2}  \langle \Delta_2 \rangle_\beta
	q_0^{(l_0,\beta)} \cdot \dots \cdot q_r^{(l_r,\beta)}
	(D,\beta) (l_r, \beta) \Delta_0
$$
The only term that could contribute a negative power of $q_r$ is 
$\beta = l_r$; however, this term is cancelled by the factor $(D, \beta) = 0$.
The terms with $q_r^0$ as factor do not contribute to the sum either since
$(l_r, \beta) = 0$. Hence we get the above formula.

\smallskip

The other cases are dealt with by very similar arguments. In the case of
$l_r \circ l_r$, the summand with $\beta = l_r$ gives $\frac{1}{q_r} l_r$
as the only term containing a negative power of $q_r$, which yields our
formula. In the product of $T_1, T_2$, the summand $\beta = l_r$
either does not appear at all if $T_1$ or $T_2$ is $\Delta_2$ because of
$k_r(l_r) = 1$; otherwise,
if both $T_1, T_2 \in H^2(V_r)$, this summand is cancelled by the factor
$(T_1, \beta)$. Now if we set $q_r =0$, in the relevant sums computing
$T_1 \circ_{V_r} T_2$ only summands
with $(\beta, l_r) = 0$ are left, i.~e. with $\beta \in H_2(V_{r-1})$.
Since the Gromov--Witten invariants of such an $\beta$ are the same, whether
they are computed in $V_r$ or in $V_{r-1}$, we get the last equation
of (3.27).

\smallskip

We have thus seen that the product extends to $F$. Since it obviously
remains associative, $F$ is an algebra over
$\Q[q_0, \dots, q_r, \allowbreak{} q_1^{-1}, \dots, q_{r-1}^{-1}]$.
So we have proved the first assertion of our lemma.

\smallskip

We want to examine the fiber $F_0:= F/q_r F$. Clearly, $q_r l_r$ is an
idempotent, so we get a splitting of $F_0$ into the image and the kernel $K$
of multiplication with $q_r l_r$. The kernel $K$ is spanned by
$H^2(V_{r-1}), \Delta_2$ and $\Delta_0 - q_r l_r$, the image is the
span of $q_r l_r$. 

\smallskip

The kernel $K$ is isomorphic to $F_0 /q_r l_r$. From the last equation in
the list (3.27) we see that $F_0 / q_r l_r$ is isomorphic to the
small quantum cohomology algebra of $V_{r-1}$, which proves the second
statement of the lemma.

\medskip

{\bf Proof of Theorem 3.6.1.} 
As mentioned already in 1.6, it is well known that the quantum
cohomology of $\P^r$ is generically semisimple. We proceed by
induction.

\smallskip

So assume we know that the quantum cohomology of $V_{r-1}$ is
generically semisimple on $H^2$. According to the lemma, this implies
that 
$$
\Spec\,F \to \allowbreak{}\Spec\, \Q[q_0, \dots, q_r,\allowbreak{}
q_1^{-1}, \dots, q_{r-1}^{-1}]
$$
is generically semisimple on the fiber $q_r = 0$. Since semisimplicity is a
Zariski--open condition, this map is generically semisimple; restricting
to $q_r \neq 0$, we see that the quantum cohomology of $V_r$ is
generically semisimple on $H^2(V_r)$.

\medskip

{\bf 3.6.3. Remark.} It is well known that there exists a full collection
of exceptional sheaves for del Pezzo surfaces, see [KuOr]. Hence the
Theorem 3.6.1 is a special case of Dubrovin's conjecture mentioned in
1.8.2.

\newpage

\centerline{\bf References}

\medskip

[Bar] S.~Barannikov. {\it Semiinfinite Hodge structures and mirror symmetry
for projective spaces.} Preprint math.AG/0010157.

\smallskip

[Beau] A.~Beauville. {\it Quantum cohomology of complete
intersections.} Mat. Fiz. Anal. Geom.  2  (1995),  no. 3--4, 384--398.
Preprint alg-geom/9501008.

\smallskip

[BrL] J.~Bryan, N.~C.~Leung. {\it The enumerative geometry of
$K3$ surfaces and modular forms.} Journ. Amer. Math. Soc.,
13:2 (2000), 371--410.

\smallskip

[Du] B.~Dubrovin. {\it Geometry and analytic theory of Frobenius
manifolds.} Proc. ICM Berlin 1998, vol. II, 315--326.

\smallskip

[DuZh] B.~Dubrovin, Y.~Zhang. {\it Bihamiltonian hierarchies
in 2D topological field theory at one--loop approximation.}
Comm.~Math.~Phys., 198 (1998), 311--361.

\smallskip

[G1] A.~Givental. {\it A mirror theorem for toric complete
intersections.} Progr. Math. 160 (1998), 141--175.

\smallskip

[G2] A.~Givental. {\it Semisimple Frobenius structures at higher genus.}
Preprint math.AG/0008067.

\smallskip

[G3] A.~Givental. {\it On the WDVV--equation in quantum $K$--theory.}
Preprint math.AG/0003158.

\smallskip

[GP] L.~G\"ottsche, R.~Pandharipande. {\it The quantum cohomology
of blow--ups of $\bold{P}^2$ and enumerative geometry.}
J.~Diff.~Geometry, 48 (1998), 61--90.

\smallskip

[Guz] D.~Guzzetti. {\it Stokes matrices and monodromy groups
of the quantum cohomology of projective space.} Comm. Math. Phys.  207 
(1999),  no. 2, 341--383.

\smallskip

[He1] C.~Hertling. {\it Multiplication on the tangent bundle.}
math.AG/9910116.

\smallskip

[He2] C.~Hertling. {\it Frobenius manifolds and moduli
spaces for hypersurface singularities.} Preprint, 2000.

\smallskip

[HeMa] C.~Hertling, Yu.~Manin. {\it Weak Frobenius manifolds.}
Int.~Math.~Res.~Notes, 6 (1999), 277--286.
math.QA/9810132.

\smallskip

[Hirsch] A.~Hirschowitz. {\it Sym\'etries des surfaces rationelles
g\'en\'eriques.} Math.~Ann. 281 (1988), 255--261.

\smallskip

[KoMa] M.~Kontsevich, Yu.~Manin. {\it
Gromov--Witten classes, quantum cohomology, and enumerative 
geometry.} Comm. Math. Phys.,
164:3 (1994), 525--562.

\smallskip

[KuOr] S.~A.~Kuleshov, D.~O.~Orlov. {\it Exceptional sheaves on Del
Pezzo surfaces.} (Russian) Izv. Ross. Akad. Nauk Ser. Mat. 58:3 (1994),
53--87; translation in Russian Acad. Sci. Izv.
Math. 44:3 (1995), 479--513. 

\smallskip

[Ma1] Yu.~Manin. {\it Three constructions of Frobenius manifolds:
a comparative study.} Asian J. Math., 3:1 (1999), 179--220
 (Atiyah's Festschrift). Preprint math.QA/9801006.

\smallskip

[Ma2] Yu.~Manin. {\it Frobenius manifolds, quantum cohomology, and moduli
spaces.} AMS Colloquium Publications, vol. 47, Providence, RI, 1999,
xiii+303 pp.

\smallskip

[Ma3] Yu.~Manin. {\it  Cubic Forms: Algebra, Geometry, Arithmetic.}
North Holland, Amsterdam, 1974 and 1986, 326 pp.

\smallskip

[MaTschi] Yu.~Manin, Yu.~Tschinkel. {\it Points of bounded height 
on del Pezzo surfaces.}
Comp. Math. 85 (1993), 315--332. Reprinted in: Yu.~Manin,
Selected Papers, World Scientific, 1996, 390--407.

\smallskip

[Sa] C.~Sabbah. {\it Hypergeometric periods for a tame polynomial.}
math.AG/9805077.

\smallskip

[TX] G.~Tian, G.~Xu. {\it On the semisimplicity of the quantum
cohomology algebras of complete intersections.} 
Math. Res. Lett.  4  (1997),  no. 4, 481--488.
Preprint alg--geom/9611035.

\enddocument